\documentclass[preprint,11pt]{elsarticle}
\biboptions{numbers,sort&compress}
\usepackage{lscape}
\usepackage{lipsum}
\usepackage{amsfonts,mathrsfs}
\usepackage{amsmath,amssymb,bbm}
\usepackage{graphicx}
\usepackage{epstopdf}
\usepackage[ruled,vlined]{algorithm2e}
\SetKwInput{KwInput}{Input}                
\SetKwInput{KwOutput}{Output}              
\usepackage{algorithmic}
\allowdisplaybreaks

\ifpdf
\DeclareGraphicsExtensions{.eps,.pdf,.png,.jpg}
\else
\DeclareGraphicsExtensions{.eps}
\fi

\newtheorem{remark}{Remark}[section]

\DeclareMathOperator{\diag}{diag}

\usepackage{float}

\usepackage{cancel}
\usepackage{subfigure}
\usepackage{bm}

\usepackage{colortbl,dcolumn}    
\usepackage{color} 
\usepackage{cancel} 
\usepackage{soul}
\usepackage{setspace}  
\usepackage{mathptmx}   
 
\newcommand{\T}{\mathsf{T}}

\newcommand{\CE}{{\mathcal{E}}}

\newcommand{\CO}{{\mathcal{O}}} 
\newcommand{\CP}{{\mathcal{P}}}

\newcommand{\IR}{{\mathbb{R}}}
\newcommand{\IC}{{\mathbb{C}}}

\usepackage[colorlinks,bookmarksopen,bookmarksnumbered,citecolor=blue,urlcolor=red]{hyperref}
\newcommand{\x}{{\bf x}}

\newcommand{\U}{{\mathbf{U}}}
\newcommand{\F}{{\mathbf{F}}}
\newcommand{\G}{{\mathbf{G}}}

\newcommand{\tol}{{\sf tol}}
\newcommand{\tn}[1]{\textnormal{#1}\ } 
\newcommand{\bmt}{\left[ \begin{array}{ccccccccccccccccccccccccccccccccccccc}}
	\newcommand{\emt}{\end{array}\right]}
\usepackage{booktabs}

\newcommand{\bmtx}{\left[ \begin{array}{c|cccccccccccccccccccccccccccccccccccc}}
	\newcommand{\emtx}{\end{array}\right]}

\newcommand{\bean}{\begin{eqnarray*}}
	\newcommand{\eean}{\end{eqnarray*}}
\newcommand{\bea}{\begin{eqnarray}}
	\newcommand{\eea}{\end{eqnarray}}
\newcommand{\eq}{\begin{equation}\begin{array}{lllllllll}}
		\newcommand{\ee}{\end{array}\end{equation}}
\newcommand{\eqn}{\begin{equation*}\begin{array}{lllllllll}}
		\newcommand{\een}{\end{array}\end{equation*}}
\usepackage{multirow,multicol}
\newcommand{\Ii}{{\mathfrak{i}}}

\renewcommand{\Re}{\mathsf{Re}}

\begin{document}

\begin{frontmatter}		
		\title{A ROM-accelerated parallel-in-time preconditioner for\\ solving all-at-once systems from evolutionary PDEs
		}
		\author[Liu]{Jun Liu}
		\ead{juliu@siue.edu} 
		\author[Wang]{Zhu Wang\corref{cor1}}
		\ead{wangzhu@math.sc.edu}
		
		\address[Liu]{Department of Mathematics and Statistics, Southern Illinois University Edwardsville, Edwardsville, IL 62026, USA.}
	\address[Wang]{Department of Mathematics, University of South Carolina, Columbia, SC 29208, USA.}
		\cortext[cor1]{Corresponding author}

\begin{abstract}
	In this paper we propose to use model reduction techniques for speeding up
		 the diagonalization-based parallel-in-time (ParaDIAG) preconditioner, for iteratively solving all-at-once systems from evolutionary PDEs. 
		 In particular, we use the reduced basis method to seek a low-dimensional approximation to the sequence of complex-shifted systems arising from Step-(b) of the ParaDIAG preconditioning procedure. 
		 Different from the standard reduced order modeling that uses the separation of offline and online stages, we have to build the reduced order model (ROM) online for the considered systems at each iteration.
		 Therefore, several heuristic acceleration techniques are introduced in the greedy basis generation algorithm, that is built upon a residual-based error indicator, to further boost up its computational efficiency.
         Several numerical experiments are conducted, which illustrate the favorable computational efficiency of our proposed ROM-accelerated ParaDIAG preconditioner, in comparison with the state of the art multigrid-based ParaDIAG preconditioner. 
\end{abstract}
        
	\begin{keyword} 
	{Parallel-in-time\sep $\alpha$-circulant preconditioner \sep model order reduction \sep reduced basis method \sep FGMRES }
	\end{keyword}

\end{frontmatter}  
		 
\section{Introduction}
With the advent of massively parallel processors, various parallel-in-time (PinT) algorithms have been developed in the last few decades for simulating time-dependent partial differential equations (PDEs) \cite{gander201550,ong2020applications}. 
Such PinT algorithms with successful implementations can provide significant speed up over the traditional sequential time-stepping schemes. 
However, the design of effective PinT algorithms is more challenging than their counterparts in space, such as spatial domain decomposition, because of the sequential nature of forward time evolution/marching. 
Thus far, there are several different types of PinT algorithms in literature\footnote{For an extensive list of PinT-related literature, see the website \url{http://parallel-in-time. org}.}, including the parareal algorithm \cite{LMT01}, the multigrid reduction in time (MGRIT) algorithm \cite{FS14}, the space-time parallel multigrid algorithm\cite{gander2016analysis},  and the more recent diagonalization-based ParaDIAG algorithm \cite{gander2020paradiag}, etc.  The mechanism of each method varies greatly, which leads to significant difference with respect to application scopes, convergence properties and parallel efficiency. 
Among them, the ParaDIAG algorithms \cite{gander2020paradiag} are built upon the  diagonalization of the time discretization matrix or its approximations within the all-at-once system arising from solving all the time steps simultaneously.
Extensive numerical results given in \cite{GW19,gander2020paradiag} indicate that the ParaDIAG algorithms have a very promising parallel efficiency for both parabolic \cite{lin2020all,wu2020diagonalization} and hyperbolic PDEs \cite{gander2019direct,wu2020parallel,WL2020SIMAX}. 
Especially, there is a sequence of sparse complex-shifted linear systems to be solved in the algorithm. These  linear systems are independent, thus can be solved by parallel computing. However, the size of the individual linear systems can be tremendous in real-world applications, as it is determined by the number of spatial degrees of freedom. Solving these systems represents the major computational cost in numerical simulations of evolutionary PDEs by the ParaDIAG algorithms. Therefore, in order to accelerate the ParaDIAG algorithms, we propose to synthesize it with model reduction techniques.

The reduced order models (ROMs) have been widely used in approximating large-scale linear and nonlinear dynamical systems, with successful applications to numerical simulations, control and optimization problems. Common model reduction techniques include but not limited to the reduced basis method (RBM), proper orthogonal decomposition (POD), dynamical model decomposition (DMD), and interpolatory methods \cite{hesthaven2015certified,HLB96,quarteroni2015reduced,antoulas2020interpolatory,antoulas2005approximation,brunton2019data}. 
Although they provide quite different ways to construct ROMs, they all are date-driven approaches and usually share a common computational strategy: (i) finding a low-dimensional approximation to the solution manifold of the underlying system's dynamics and building a ROM at an offline stage; (ii) using the ROM for fast simulations at an online stage. The efficacy of the low-dimensional approximation can be quantified by the Kolmogorov $n$-width \cite{pinkus2012n,binev2011convergence}. 
Successful ROMs can dramatically reduce the online simulation cost.  
Recently, the ROMs have also been integrated with preconditioners in Krylov-subspace iterative solvers. For instance, preconditioned Conjugate Gradient (PCG) methods are developed using a POD-based preconditioner in  \cite{markovinovic2006accelerating,pasetto2017reduced} and a RBM-based one in \cite{nguyen2018reduced}, which are able to speed up the traditional iterative methods in serial computing. 

In this work, we pursue in this direction and improve the efficiency of the ParaDIAG algorithms by online reduced order modeling. 
We focus on addressing the sequence of linear systems in Step-(b) of ParaDIAG via building their low-dimensional approximation by the RBM. In particular, a greedy algorithm is used to find a reduced basis and the Galerkin projection is applied in constructing the ROM. Once the ROM is constructed, the computational complexity for solving the systems only depends on the dimension of the ROM, which is much smaller than the number of degrees of freedom in space and time of the full order models. 
Since the build stage is part of the approach, we propose several algorithmic improvements on the greedy basis generation. 
The proposed approach would lead to significant computational savings in Step-(b) and, as a consequence, notably decreases the overall computational cost of the ParaDIAG algorithms.   

The rest of the paper is organized as follows.
In the next section, a block $\alpha$-circulant type ParaDIAG algorithm is briefly reviewed for solving an unsteady convection-diffusion equation based on a standard upwind finite difference scheme.
Our proposed ROM-based ParaDIAG preconditioner for approximately solving the linear systems in Step-(b) is introduced in Section \ref{ROMsolver}, where a greedy algorithm for selecting the reduced basis is explained in detail and brief discussions on computational complexity are presented. Extensive 1D and 2D numerical examples are presented in Section \ref{secNum} to demonstrate the promising performance of our proposed ROM-based ParaDIAG preconditioner in contrast with the multigrid-based ParaDIAG preconditioner.
Finally, some conclusions are drawn in the last section.
\section{The ParaDIAG algorithm for fully discretized evolutionary PDEs}
\label{ModelFDM}
In this work, we consider the unsteady convection-diffusion equation \cite{morton1996numerical}:
\begin{equation}\label{modelPDE}
	\begin{cases}
		u_{t}(\x,t)= \nabla\cdot(a(\x)\nabla u)-\bm c(\x)\cdot  \nabla u(\x,t) +f(\x,t), &\tn{in} \Omega\times(0,T],\\
		u(\x,t)=g(\x,t),&\tn{on} \partial\Omega\times(0, T],\\
		u(\x,0)=u_0(\x),&\tn{in} \Omega\times\{0\}, 
	\end{cases}
\end{equation} 
where $\Omega\subset \IR^d$ is the spatial domain, $a(\x)\ge a_0>0$ is a variable diffusion coefficient, $\bm c(\x)$ is a convection or velocity field, $f$ is a source term, $g$ is a Dirichlet boundary condition, and $u_0$ is an initial condition. 
To illustrate the ParaDIAG algorithm, we restrict ourselves to the case in which $\Omega=[0,1]^2$ $(d=2)$ and  discretize the space and time using finite difference methods (FDM). However, we emphasize the proposed algorithm can be applied to problems defined in a one-dimensional or three-dimensional domain; it can also be applied to problems defined in an irregular domain for which finite element method can be employed for spatial discretization. 
Furthermore, 
we assume all data are regular enough such that the solution is unique and sufficiently smooth to assure the convergence of finite difference schemes. We refer to the monographs \cite{morton1996numerical,roos2008robust,hackbusch2017elliptic} for more discussion on various different discretization schemes.
For brevity, we adopt the commonly used notations to denote the following reshaping operations: matrix-to-vector $\texttt{vec}$ and vector-to-matrix $\texttt{mat}$, as
defined in \cite[p. 28]{golub2012matrix}.
We will also use the well-known matrix Kronecker product property that
$\texttt{vec} (AXB)= (B^\T\otimes A) \texttt{vec}(X)$ for any matrices $A,B,X$ of compatible sizes.
 
For the discretization of \eqref{modelPDE}, we use the backward Euler scheme in time,
the conservative central finite difference scheme for the diffusion term,  and the upwind scheme for the convection term  \cite{hackbusch2017elliptic,LeVeque2007}. 
Given a positive integer $N$, let $h=1/(N+1)$ be the spatial mesh size and the 2D square domain $\overline\Omega=[0,1]^2$ be partitioned uniformly by grid points 
$(x_i=ih, y_j=jh)$ with $0\leq i, j \leq N+1$.  
We denote the set of all interior grid points by $\Omega_h=\{(x_i,y_j)\}_{i,j=1}^{i,j=N}$ and also define the cell center points $(x_{i\pm 1/2}=(i\pm 1/2)h, y_{j\pm 1/2}=(j\pm 1/2)h)$ with $1\leq i, j \leq N$.
Given the final time $T>0$ and another positive integer $K$, let $\tau=T/K$ be the time step size and the time domain $[0,T]$ be uniformly divided by grid points 
$t_k=k\tau$ with $0\leq k \leq K$. 
Let $U_{i,j}^k$ represent the finite difference approximation of $u(x_i,y_j,t_k)$. 
Assuming $\bm c(\x)=(p(x,y), q(x,y))$, 
the full finite difference discretization of the PDE (\ref{modelPDE}) on
the set $\Omega_h$ of interior grid points  reads: for $1\le i,j\le N$, $1\le k\le K$,  
\begin{align} \label{FDM}
\frac{U_{i,j}^{k}-U_{i,j}^{k-1}}{\tau}=&
\frac{1}{h}\left(a_{i+1/2,j} \frac{U_{i+1,j}^{k}-U_{i,j}^{k}}{h}-a_{i-1/2,j} \frac{U_{i,j}^{k}-U_{i-1,j}^{k}}{h}\right) \nonumber\\ &+\frac{1}{h}\left(a_{i,j+1/2} \frac{U_{i,j+1}^{k}-U_{i,j}^{k}}{h}-a_{i,j-1/2} \frac{U_{i,j}^{k}-U_{i,j-1}^{k}}{h}\right) \nonumber\\
&-\left(p_{i,j}^+\frac{U_{i,j}^{k}-U_{i-1,j}^{k}}{h}
+p_{i,j}^-\frac{U_{i+1,j}^{k}-U_{i,j}^{k}}{h} \right)\nonumber\\ &-\left(q_{i,j}^+\frac{U_{i,j}^{k}-U_{i,j-1}^{k}}{h}
+q_{i,j}^-\frac{U_{i,j+1}^{k}-U_{i,j}^{k}}{h} \right) 
 +f_{i,j}^{k}, 
\end{align} 
where $a_{i\pm 1/2,j\pm 1/2}=a(x_{i\pm 1/2}$, $y_{j\pm 1/2})$,
 $p_{i,j}=p(x_i,y_j)$, $q_{i,j}=q(x_i,y_j)$,
and $f_{i,j}^{k}=f(x_i,y_j,t_{k})$. 
Here we defined the short notations $z^+:=\max\{z,0\}$ and $z^-:=\min\{z,0\}$
for $p_{i,j}$ and $q_{i,j}$ to formulate the upwind scheme. 
Under suitable regularity assumptions as discussed in \cite[Thm. 5.17]{hackbusch2017elliptic}, the above full scheme (\ref{FDM}) has a first-order accuracy in both space and time.

Let $\U^k,\F^k$ be column vectors containing the ordered values of $U_{i,j}^{k},f_{i,j}^{k}$ over all the interior grids in $\Omega_h$ at time $t_k$, respectively. 
Upon enforcing
the Dirichlet boundary condition by shifting boundary nodes to the right-hand-side,
the above finite difference scheme (\ref{FDM}) can be written into 
the following sequential time stepping formulation (marching from the given $\U^0$ to $\U^1$, then to $\U^2$, till the last time step $\U^K$): 
\begin{align} \label{BEtimestepping}
	(\tau^{-1} I_h- L_h) \U^k= \tau^{-1}\U^{k-1}+ \F^k+\G^k, \quad 1\le k\le K,
\end{align} 
where $I_h\in\IR^{N^2\times N^2}$ is an identity matrix, $L_h$ is the coefficient matrix representing the discretized spatial differential operators (involving $a,p,q$), $\U^0$ is given by the initial condition $u_0(x,y)$, and $\G^k$ comes from the Dirichlet boundary conditions.
Explicitly, the coefficient matrix $L_h$ has the following  expression in terms of Kronecker products:
\begin{align} \label{Ahmat}
L_h=&-\frac{1}{h^2}\left[ (I\otimes D^\T )A_x (I \otimes D)+(D^\T\otimes I)A_y (D \otimes I)\right] \nonumber\\
&-\frac{1}{h} \left[ P^+ (I\otimes D^+) +P^- (I\otimes D^-)+
 Q^+ ( D^+\otimes I) +Q^- (D^- \otimes I)
\right],
\end{align}	
where $I\in\IR^{N\times N}$ is an identity matrix,
\begin{align*}
	D&= \begin{bmatrix}
		1 & & & &\\
		-1 &1  & & &\\
		0 &-1 &1 & &\\
		&\ddots &\ddots  &\ddots &\\
		& &0  &-1 &1\\ 
			& &  &0 &-1\\ 
	\end{bmatrix}\in\IR^{(N+1)\times N},\quad D^+=D(1:N,:),\quad D^-=D(2:N+1,:),\\
A_x&=\diag\left(\texttt{vec}\left([a_{i-1/2,j}]_{i=1,j=1}^{N+1,N}\right)\right),\quad
A_y=\diag\left(\texttt{vec}\left([a_{i,j-1/2}]_{i=1,j=1}^{N,N+1}\right)\right),\\
P^{\pm}&=\diag\left(\texttt{vec}\left([p^{\pm}_{i,j}]_{i=1,j=1}^{N,N}\right)\right),\quad 
Q^{\pm}=\diag\left(\texttt{vec}\left([q^{\pm}_{i,j}]_{i=1,j=1}^{N,N}\right)\right).
\end{align*}
Here $Z=[z_{i,j}]_{i=1,j=1}^{m,n}$ defines an $m\times n$ matrix with $z_{i,j}$ as its $(i,j)$-th entry,
$\texttt{vec}(Z)$ reshapes an $m\times n$ matrix $Z$ into a
column vector $\bf Z\in \IR^{mn}$ by stacking all columns together one after another, and
$\diag(\bf Z)$ denotes a diagonal matrix with the vector $\bf Z$ as its main diagonal entries. Notice that $D^+$ and $D^-$ are defined by selecting the first and last $N$ rows of
the $(N+1)\times N$ rectangular difference matrix $D$, respectively.

Different from the sequential time stepping method marching from the initial time step to the final,
the all-at-once approach tries to simultaneously solve all time steps in one-shot without explicit time marching. More specifically, by stacking all the $K$ systems in (\ref{BEtimestepping}) together and using Kronecker product notation, we can obtain  the following all-at-once non-symmetric $N^2K\times N^2K$ linear system  
\begin{align}
	 (\tau^{-1} B\otimes I_h+I_t\otimes (- L_h) ) \U= \F,
	\label{eq:bigsys}
\end{align}
where  $I_t\in\IR^{K\times K}$ is an identity matrix,
\begin{align}
	B= \begin{bmatrix}
		1 & & & &\\
		-1 &1  & & &\\
		0 &-1 &1 & &\\
		&\ddots &\ddots  &\ddots &\\
		& &0  &-1 &1\\ 
	\end{bmatrix}\in \IR^{K\times K},\quad 
\U=\bmt \U^1\\ \U^2\\ \U^3\\ \vdots \\ \U^K \emt,\quad 
\F=\bmt  (\F^1+\G^1)+\frac{1}{\tau}\U^0 \\  (\F^2+\G^2)\\  (\F^3+\G^3)\\ \vdots \\  (\F^K+\G^K) \emt.
\end{align} 
By exploiting the Toeplitz matrix $B$ due to the uniform time step size, 
the Strang-type block $\alpha$-circulant preconditioner \cite{McDonald_2018,lin2020all} for (\ref{eq:bigsys}) has the following   Kronecker product form
\begin{align}
	\CP_{\alpha}:= \tau^{-1} C_{\alpha}\otimes I_h+I_t\otimes (-  L_h),    
\end{align}
where the Toeplitz matrix $B$ is replaced by a Strang-type $\alpha$-circulant matrix (with $\alpha\in (0,1)$)
\begin{align}
	C_{\alpha}= \begin{bmatrix}
		1 & & & &-\alpha\\
		-1 &1  & & &\\
		0 &-1 &1 & &\\
		&\ddots &\ddots  &\ddots &\\
		& &0  &-1 &1\\ 
	\end{bmatrix} \in\IR^{K\times K}.
\end{align}
It is shown in \cite{lin2020all} that, for a sufficiently small $\alpha=O(\sqrt{\tau})\in (0,1)$, the GMRES equipped with the above preconditioner $\CP_{\alpha}$ can solve \eqref{eq:bigsys} effectively and it achieves  a provable mesh-independent convergence rate. 
Here $C_{\alpha}$ is indeed a rank-one perturbation of $B$
satisfying $\lim_{\alpha\to 0}\|C_{\alpha}-B\|=0$,
which implies the preconditioner $\CP_{\alpha}$ is expected to lead faster convergence rate as $\alpha$ gets smaller if not considering the possible round-off errors due to inverting $\CP_{\alpha}$ during the preconditioning step.
In practice,  a small $\alpha=0.01$ seems to be effective.

Let $\mathbb{F}=\frac{1}{\sqrt{K}}\left[\omega^{(l_1-1)(l_2-1)}\right]_{l_1, l_2=1}^{K}$ (with $\Ii=\sqrt{-1}$ and $\omega=e^{\frac{2\pi{\Ii}}{K}}$) 
be the discrete Fourier matrix and  define a diagonal matrix
$
	\Gamma_\alpha=\mathrm{diag}\left(1,\alpha^{\frac{1}{K}},\cdots,\alpha^{\frac{K-1}{K}}\right)
$ for $\alpha\in (0,1)$.
The above so-called $\alpha$-circulant matrix $C_{\alpha}$   can be explicitly  diagonalized \cite{Bini_2005} according to its spectral decomposition
\begin{equation}
C_{\alpha}=V D V^{-1}, 
\end{equation}
where   $V=\Gamma_\alpha^{-1}\mathbb{F}^*$ and 
$D =\text{diag}(d_{1},\cdots,d_n,\cdots,d_{K}):=\mathrm{diag}\left(\sqrt{K}\mathbb{F}\Gamma_\alpha C_{\alpha}(:,1)\right)
$ 
with $C_{\alpha}(:,1)$ being the first column of $C_{\alpha}$. More specifically, let $\theta_n=2(n-1)\pi/K\in [0,2\pi)$ with $1\le n\le K$, the explicit expressions for the $K$ eigenvalues of $C_\alpha$ are
\begin{equation} \label{eigCalpha}
d_n= 1-\alpha^{\frac{1}{K}} e^{\Ii \theta_n}= (1-\alpha^{\frac{1}{K}}\cos\theta_n)-\Ii  \alpha^{\frac{1}{K}}\sin\theta_n,\quad 1\le n\le K.
\end{equation}
With the above explicit diagonalization $C_{\alpha}=V D V^{-1} $, we can easily factorize $\CP_\alpha$ in into
\[
 {\CP}_\alpha= \underbrace{(V\otimes I_h)}_{\tn{Step-(a)}}\
 \underbrace{\left({\tau}^{-1} D\otimes I_h+I_t\otimes (-  L_h)\right)}_{\tn{Step-(b)}}\ 
 \underbrace{(V^{-1}\otimes I_h)}_{\tn{Step-(c)}},
\]
where the matrix of Step-(b) has the following block diagonal structure
\[{\tau}^{-1} D\otimes I_h+I_t\otimes (-  L_h)= 
\begin{bmatrix}
	(d_{1}/\tau)I_h-  L_h & & & & \\
	  &\ddots  & & &\\
	  &  &(d_{n}/\tau)I_h- L_h & &\\
	&  &  &\ddots  &\\
	& &   &  &(d_{K}/\tau)I_h-  L_h\\ 
\end{bmatrix}.
\] 
Such a block diagonal structure is the key for designing the diagonalization-based ParaDIAG algorithms.  
Suppose the preconditioner $\CP_{\alpha}$ is used in a preconditioned Krylov subspace method, it is required to compute or approximate the preconditioning step $\CP_\alpha ^{-1} s$ at each iteration, where $s\in\IR^{N^2 K}$ is the residual vector.
Let $S:=\texttt{mat}(s)$,  the preconditioning step $z:=\CP_{\alpha}^{-1}s$ can be implemented via the following 3 steps:
\begin{equation}\label{3step}
	\begin{split}
		&\text{Step-(a)}\quad ~~S_1=S(V^{-1})^{\T}  ,\\
		&\text{Step-(b)}\quad ~~S_{2}(:,n)=((d_{n}/\tau)I_h- L_h)^{-1}S_{1}(:,n),\quad ~n=1,2,\dots,K,\\
		&\text{Step-(c)}\quad ~~z=\texttt{vec}(S_2V^\T) ,\\
	\end{split}
\end{equation} 
where 
$S_1(:,n)$ and $S_2(:,n)$ denote the $n$-th column of $S_1$ and $S_2$,  respectively. Here we used $\texttt{vec} (AXB)= (B^\T\otimes A) \texttt{vec}(X)$.

Since $V=\Gamma_\alpha^{-1}\mathbb{F}^*$ and $V^{-1}=\mathbb{F}\Gamma_\alpha$,  Steps-(a) and (c) in (\ref{3step}) can be computed efficiently via FFT in time direction at the computational complexity $\CO(N^2 K\log K)$ operations. Step-(b) needs to solve $K$ complex-shifted elliptic systems of size $N^2\times N^2$ with different right-hand-sides. 
Since these systems are independent of each other, parallel computing can be used to significantly reduce the wall-clock time. 
If $N$ is large, these systems can be costly if solved by direct sparse solvers. 
Therefore, another efficient iterative solver has to be applied in Step-(b). 
This essentially leads to inexact GMRES \cite{simoncini2003theory} or flexible GMRES (FGMRES) \cite{saad1993flexible,simoncini2002flexible}.
For instance, the authors in \cite{lin2020all} employ one geometric multigrid V-cycle with ILU smoother  to approximately solve the sparse linear systems in Step-(b), which needs slightly more outer GMRES iterations than solving these systems by the costly sparse direct solver.
If the preconditioned GMRES converges in $l$ iterations for a given tolerance, the overall computational cost of such a multigrid-based ParaDIAG (ParaDIAG-MG) preconditioner for solving \eqref{eq:bigsys} is about $\CO(N^2 K\log K+ lN^2K)$ operations.
Other efficient iterative solvers, 
such as domain decomposition algorithms \cite{leng2019additive,gander2019class,taus2020sweeps}, can also be used                                                                                                                                                                                                                                                                                                                                                                     in Step-(b), but the overall computational complexity remains high. 
Thus, it is natural to consider model reduction techniques to derive a low-dimensional surrogate model for these systems that could greatly improve the computational efficiency. 
\subsection{A motivating example}\label{sec:me1} 

\begin{figure}[!htb]
	\begin{center}
		\includegraphics[width=1\textwidth]{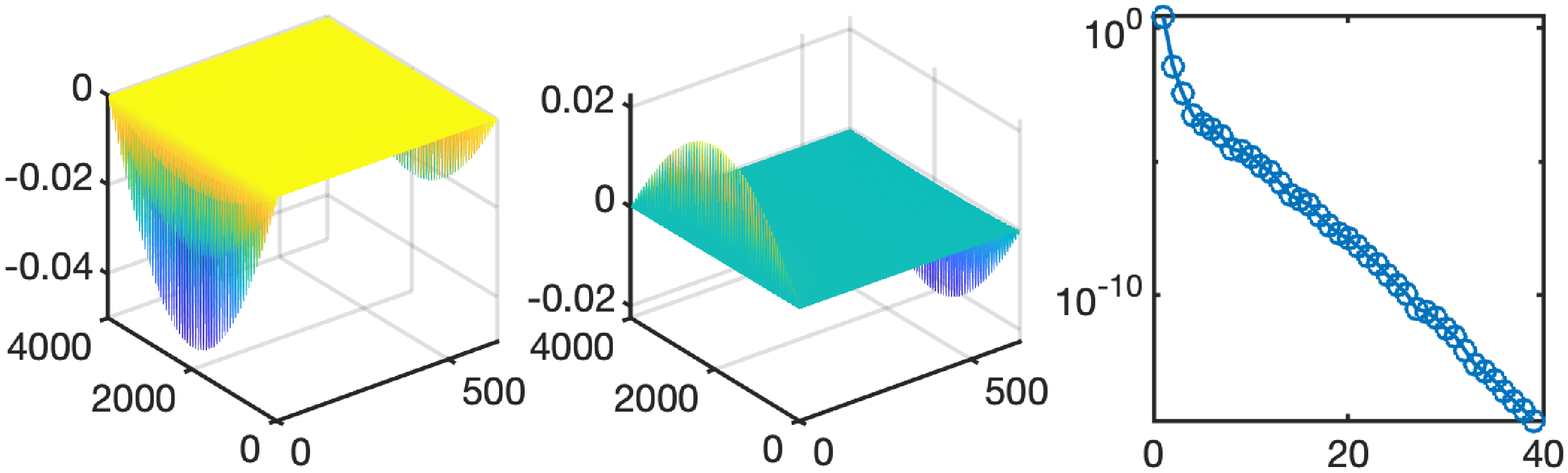} 
		\includegraphics[width=1\textwidth]{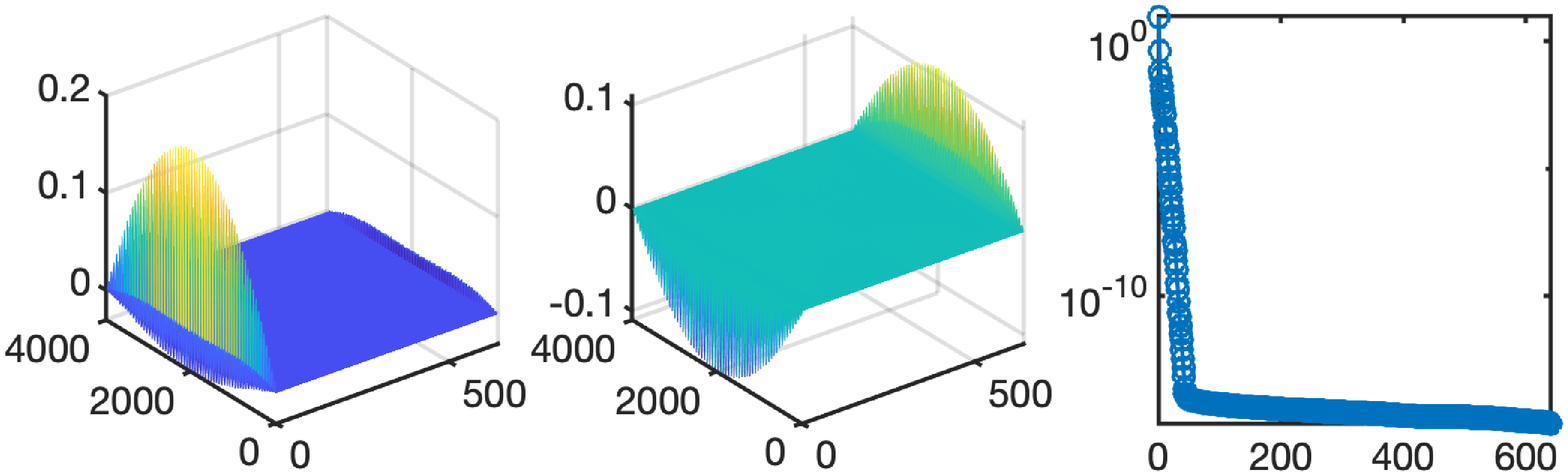} 
		\includegraphics[width=1\textwidth]{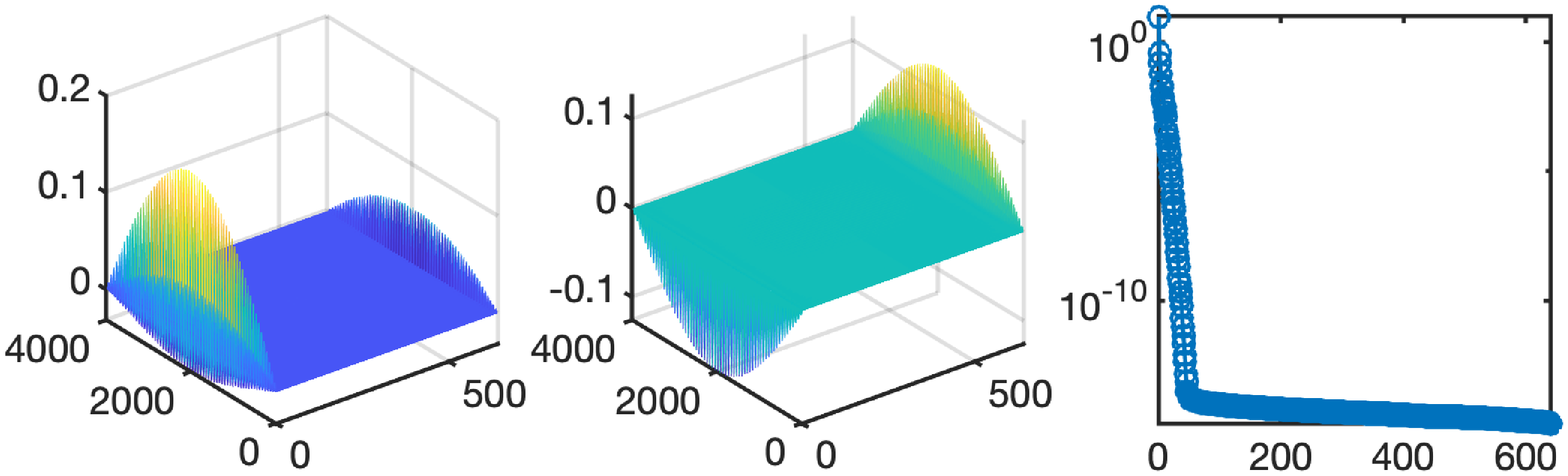} 
		\includegraphics[width=1\textwidth]{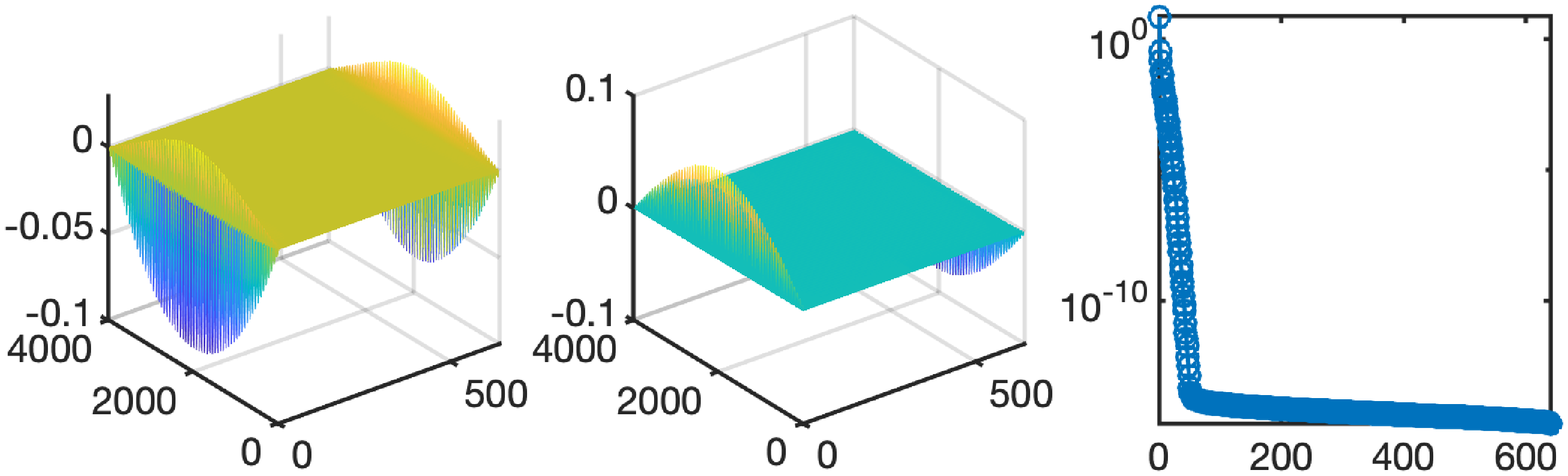} 
	\end{center}
	\caption{The visualization of Step-(b) outputs $S_2$ at iterations 1 to 4 (from top to bottom): (left column) real part of $S_2$, (middle column) imaginary part of $S_2$, (right column) singular values of $S_2$ greater than $10^{-15}$.
	}
	\label{plot_2D_ex1_S2}
\end{figure}

Before presenting the proposed ROM-based solver for Step-(b), we motivate it by considering the heat equation with a constant diffusion coefficient, defined on the unit square domain. In particular, the equation is obtained from \eqref{modelPDE} by setting $a(x,y)=0.01$, $\bm c(x,y)= (0,0)$, $f(x,y,t)=0$, $g(x,y,t)=0$ and $u_0(x,y)=x(x-1)y(y-1)$: 
  \begin{equation}\label{heatPDE}
  	\begin{cases}
  		u_t=0.01 \Delta u, &\tn{in} \Omega\times(0,T),\\
  		u=0,&\tn{on} \partial\Omega\times(0, T), \\
		u=x(x-1)y(y-1),&\tn{in} \Omega \times \{0\},
  	\end{cases}
  \end{equation}  
with the final time $T=10$.
The known exact solution is used for measuring approximation errors. 
We set $h= \tau= 1/64$ in the FDM discretization and use one multigrid V-cycle to approximately solve the linear systems of Step-(b) in serial.

The FGMRES preconditioned by the ParaDIAG-MG with ILU smoothing converges after 4  iterations, which achieves the accuracy $5.7\times 10^{-5}$. Meanwhile, Step-(b) dominates the total computational cost: the simulation is completed in 6.2~seconds, while 5.5~s of them is spent in solving Step-(b) that takes over 88\% of the total CPU time. 
Hence the central task for improving the overall efficiency is to speed up the system solving in Step-(b).

The Step-(b) outputs $S_2$ during the four iterations are shown in Fig. \ref{plot_2D_ex1_S2}: real part in the left column and imaginary part in the middle column, where relatively large structures are observed mainly near both ends. We perform the singular value decomposition (SVD) of $S_2$ and display the singular values of $S_2$ greater than $10^{-15}$ in the right column of Fig. \ref{plot_2D_ex1_S2}. It is seen that the singular values decay quickly, which indicates the solution manifolds have low-dimensional structures. This further motivates us to use a model reduction technique such as RBM to extract the low-rank trail spaces and accelerate the computation of Step-(b) by solving reduced order systems. 
We remark that the exponential or algebraic decay rate of singular values of $S_2$
highly depends on the underlying physical model,
which usually shows much slower decay rate in convection-dominated problems.

\section{ROM-based Fast Solver for Step-(b)}
\label{ROMsolver}
The independent linear systems to be solved in Step-(b) can be regarded as discrete systems  associated to a parametric PDE (pPDE) with a varying source term. Solving them is equivalent to querying numerical solutions of the pPDE for multiple times. Reduced basis method (RBM), as one of the popular model reduction techniques, is designed for efficiently completing such a computation task \cite{hesthaven2015certified,quarteroni2015reduced,barrault2004empirical,grepl2007efficient}. Typically, it has two separate computation stages: (i) at the offline stage, a low-dimensional subspace is determined that could approximate the pPDE solution manifold subject to a user-defined error tolerance, and  a reduced order model is built on this subspace; (ii) then at the online stage, the ROM as a surrogate model to the original pPDE is simulated at a very cheap cost. However, in the context of our considered ParaDIAG preconditioning, it is not straightforward to construct the ROM offline, because it is difficult to parametrize the varying source terms. Moreover, the source terms consist of the residual vectors that would change very irregularly during the GMRES iterations. As a result, we have to move the offline stage of building the ROM to the online stage. However, if the solution stays within a low dimensional manifold, one can still expect the total computational cost is  significantly less than the full order simulations,
but some special attention are needed to reduce the computation cost of the `offline' stage during each iteration.

We first recast the systems in Step-(b) into a general form: for $n = 1, \ldots, K$,
\begin{equation}
A_n x_n:=((d_n/\tau) I_h - L_h) x_n = b_n,
\label{eq:pPDEs1}
\end{equation}
where $A_n\in\IC^{N^2\times N^2}$ with $d_n$ being a complex number. 
The model is referred to as the full order model (FOM). 
They have a continuous pPDE counterpart: 
\begin{equation}
\left(d(\theta)/\tau - L\right) v(\x; \theta) = b(\x; \theta),\quad \tn{for} \theta\in [0,2\pi),
\label{eq:pPDE}
\end{equation}
where 
$d(\theta)= 1-\alpha^{\frac{1}{K}} e^{\Ii \theta}$, 
$L(u)=\nabla\cdot(a(\x)\nabla u)$ if there is no convection term in \eqref{modelPDE}; and  
$L(u)=\nabla\cdot(a(\x)\nabla u)-\bm c(\x)\cdot  \nabla u(\x,t) +\frac{h}{2} [p(x,y) u_{xx} + q(x,y) u_{yy}]$ if there is convection. 
Note that the extra dissipation appears in the latter case due to the use of upwind finite difference \cite{quarteroni2008numerical}. 
However, we shall focus on the discrete set of equations \eqref{eq:pPDEs1}, corresponding to the pPDE at the prescribed discrete parameter set $\{\theta_1,\ldots, \theta_K\}$ specified in \eqref{eigCalpha}.  
For a long time simulation, we expect $K$ to be very large and hence give rise to a large discrete parameter set.

Assume there exists a reduced space spanned by the basis vectors $\{\phi_1, \phi_2, \ldots, \phi_r\}$
and let $\Phi_r = [\phi_1, \phi_2, \ldots, \phi_r]\in\IC^{N^2\times r}$, the reduced state can be defined by $\widetilde{x}_n = \Phi_r \widehat{x}_n$. 
Replacing the state in \eqref{eq:pPDEs1}  with the reduced approximation $\widetilde{x}_n\approx x_n$ and applying the Galerkin projection, we obtain the following reduced system: 
to find $\widehat{x}_n$ such that 
\begin{equation}
\Phi_r^* A_n  \Phi_r \widehat{x}_n = \Phi_r^* b_n, 
\label{eq:pPDEs2}
\end{equation}
where $\Phi_r^*$ is the conjugate transpose of $\Phi_r$. 
Once $\widehat{x}_n$ is found, the original state variable can be approximated by the reduced state.  
The size of this reduced system matrix $\Phi_r^* A_n  \Phi_r$ is  $r\times r$. As $r$ is much less than $N^2$, it can be solved cheaply using a direct solver with the complexity of $\CO(r^3)$.

The remaining task is to efficiently determine the reduced basis matrix $\Phi_r$. For this purpose, RBM uses a greedy strategy based on {\em a posteriori} error estimations. Such error estimations are usually derived from the equivalence between the error norm and a dual norm of the residual vector. 
Next, we focus on the discrete systems in Step-(b) and derive an error indicator by establishing the relation between the unknown error norms and measurable residual norms. 
 
\section*{A residual-based error indicator}

	Note that $\alpha\in(0,1)$, the constant $d_n= 1-\alpha^{\frac{1}{K}} e^{\Ii \theta_n}$ given in (\ref{eigCalpha}) has a positive real part 
	\begin{equation}
		\Re( d_n)= (1-\alpha^{\frac{1}{K}}\cos\theta_n)\ge  (1-\alpha^{\frac{1}{K}}) > 0,\label{eq:eig}
	\end{equation}
which implies the matrix $A_n=(d_n/\tau)I_h - L_h$ 
is strictly diagonally dominant by rows and columns
since $(-L_h)$ is irreducibly diagonally dominant \cite[Thms. 5.17]{hackbusch2017elliptic}.
Define two positive constants
\begin{align}
\beta_{R}: 
&=\min_{1\le k\le N^2} \left(|(d_n/\tau)- (L_h)_{k,k}|-\sum_{j\ne k} |(L_h)_{k,j}| \right)
\end{align}
and 
\begin{align}
\beta_C:=\min_{1\le k\le N^2} \left(|(d_n/\tau)- (L_h)_{k,k}|-\sum_{j\ne k} |(L_h)_{j,k}| \right).
\end{align}
Since $\Re( d_n)>0$ and $-(L_h)_{k,k}>0$, there obviously holds
$$|(d_n/\tau)- (L_h)_{k,k}|\ge \Re(d_n/\tau)+|(L_h)_{k,k}|,$$
 and $|(L_h)_{k,k}|\ge\sum_{j\ne k} |(L_h)_{k,j}|$ due to the row weakly diagonal dominance of $(-L_h)$, which leads to
\begin{align}
	\beta_{R}
	\ge \min_{k} \left(\Re(d_n/\tau)+|(L_h)_{k,k}|-\sum_{j\ne k} |(L_h)_{k,j}| \right)
	\ge \Re( d_n)/\tau. 
\end{align} 
Similarly there holds $\beta_C\ge \Re( d_n)/\tau $ because of the column weakly diagonal dominance of $(-L_h)$.
By a lower bound of the smallest singular value \cite[Cor. 2]{varah1975lower}, there holds 
 \begin{equation} \label{Anbound}
 	\|A_n^{-1}\|_2\le \frac{1}{\sqrt{\beta_R \beta_C}}
 	\le \frac{\tau}{\Re( d_n)}
 	\le \frac{\tau}{1-\alpha^{{1}/{K}}}
 	=  \frac{\tau}{1-\alpha^{{\tau}/{T}}}\le \frac{T}{1-\alpha},
 	\end{equation}	
 where we used the fact that $\phi(\tau):=\frac{\tau}{1-\alpha^{{\tau}/{T}}}$ is an increasing function of $\tau\in (0,T]$. In fact, by the inequality $\ln(1+z)<z$ for $z>0$, a simple calculation
 can  verify $\phi'(\tau)>0$, that is
 \[ 
 \phi'(\tau)
 =\frac{ 1-\alpha^{{\tau}/{T}}(1+\ln(\alpha^{-{\tau}/{T}}))}{ (1-\alpha^{{\tau}/{T}})^2}
 >\frac{ 1-\alpha^{{\tau}/{T}}(1+(\alpha^{-{\tau}/{T}}-1))}{ (1-\alpha^{{\tau}/{T}})^2}
 =0.
 \]
 The bound (\ref{Anbound}) shows that some $A_n$ may become more ill-conditioned (or as ill-conditioned as $L_h$) as the chosen fixed parameter $\alpha$ gets closer to $1$.
 The second inequality also explains why the linear systems with $\Re(d_n)\approx 0$ (or $\cos(\theta_n)\approx 1$) are more difficult to approximate since their condition numbers are larger.
 We mention that the limiting case $\alpha=1$ may lead to much slower convergence rate of preconditioned GMRES (see Example 2 in \cite{lin2020all}) and hence not further considered here (see below Remark 1 for discussion of the pure diffusion situation with an upper bound covers the case $\alpha=1$).

Denote the reduced approximation error by $e_n = x_n-\widetilde{x}_n$, and define the residual vector by $r_n = b_n - A_n \widetilde{x}_n$. 
The error representation reads: 
\begin{equation}
A_n e_n =A_n(x_n-\widetilde{x}_n)=A_n x_n-A_n \widetilde{x}_n= b_n - A_n \widetilde{x}_n = r_n. 
\label{eq:err1}
\end{equation}
which implies 
\begin{equation}
\|e_n\|_2=\|A_n^{-1}r_n\|_2\le \|A_n^{-1}\|_2\|r_n\|_2\le
 \frac{\tau}{\Re( d_n)} \|r_n\|_2 \leq \frac{T}{1-\alpha}\|r_n\|_2. 
\label{eq:errb}
\end{equation}
The right-hand-side term provides an upper bound of the error norm
and the denominator $(1-\alpha)$ suggests that choosing a smaller $\alpha$ close to $0$ leads to tighter estimate. 
Numerically, we find that $\alpha=10^{-2}$ works very well for all tested examples.
Although it is possible to directly take the above upper bounds (e.g. $ \frac{\tau}{\Re( d_n)} \|r_n\|_2$) as an absolute error indicator, to be compatible with the used outer FGMRES stopping criterion, we will use the following relative residual norm 
\begin{equation}
	\CE_n:=\frac{\|r_n\|_2}{\|b_n\|_2}
	\label{errind}
\end{equation}
as an error indicator for estimating the errors and guiding the greedy  algorithm for searching reduced basis vectors. 
During the process, we also obtained the approximate solutions to Step-(b) simultaneously.
 
\begin{remark}
	The above uniform upper bound (\ref{Anbound}) may be improvable, since it does not explicitly depend on the ellipticity constant $a_0>0$ yet.
	For example, consider the pure diffusion equation with $\bm c(\x)=(0,0)$
	and constant coefficient $a(x,y)=a_0>0$, $(-L_h)$ is symmetric positive definite and there holds \cite[Thm 4.34]{hackbusch2017elliptic}
	\[
	\lambda_{\min}(-L_h)=\|(-L_h)^{-1}\|^{-1}_2\ge 16a_0,
	\]
	where $\lambda_{\min}(-L_h)$ denotes the smallest eigenvalue of $(-L_h)$. 
	Then from  $A_n e_n=r_n$ with $A_n=(d_n/\tau)I_h - L_h$, we can get 
	\[
	(d_n/\tau) e_n^*e_n+e_n^*(-L_h)e_n=e_n^* A_n e_n=e_n^* r_n,
	\]
	and since $\Re(d_n)>0$ and $e_n^*(-L_h)e_n\ge \lambda_{\min}(-L_h)\|e_n\|_2^2\ge 16a_0\|e_n\|_2^2>0$ there holds
	\[
  (\Re(d_n/\tau) +16a_0)\|e_n\|_2^2
  \le \Re(d_n/\tau) e_n^*e_n+e_n^*(-L_h)e_n
  \le |e_n^*A_n e_n|=|e_n^*r_n|\le \|e_n\|_2  \|r_n\|_2.
	\]
	This leads to the following slightly improved error estimate
	\[
	 \|e_n\|_2\le \frac{1}{\Re(d_n/\tau) +16a_0} \|r_n\|_2\le \frac{\tau}{(1-\alpha^{\tau/T})+16a_0 \tau}  \|r_n\|_2
	  \le \frac{T}{(1-\alpha)+16a_0 T}\|r_n\|_2\le \frac{1}{16a_0}\|r_n\|_2,
	\]
	which is tighter than the bound (\ref{eq:errb}) when $a_0\gg 0$ and $T$ is large and is also applicable to the limiting case $\alpha=1$.
	For a general convection-diffusion equation, the spatial discretization matrix $(-L_h)$ is not symmetric any more and the above Rayleigh quotient-based arguments utilizing a known lower bound of the smallest eigenvalue would not be valid in general. 
\end{remark}
\section*{A greedy algorithm for basis generation and its complexity} 
The reduced basis generation strategy is summarized in Algorithm \ref{alg:basis}. 
At each iteration, the algorithm identifies a problem that would yield the worst {\em a posteriori} error, and generates a new basis from the associated FOM solution. As the reduced space gets richer, the reduced approximation would become more accurate. When the accuracy meets the user-defined tolerance ${\sf tol_{ROM}}$, the process would terminate and return the approximation solutions to \eqref{eq:pPDEs1}.
Note that Step (b) is performed inside the FGMRES iteration, for which we only need to provide a reasonable approximate solution. In our implementation, ${\sf tol_{ROM}}$ is set to be square root of the tolerance of FGMRES iterations. We refer to \cite{simoncini2002flexible,Giraud2007} for further analysis on how the inexact solving of preconditioned systems affects the convergence of FGMRES or relaxed GMRES.

\begin{algorithm}[!htb]
\SetAlgoLined
\KwInput{ ${\sf \tol_{ROM}}$, $K$ and $r=1$}
\KwOutput{ reduced basis matrix $\Phi_r$ and solutions to \eqref{eq:pPDEs1} $\{{x}_1, \ldots, {x}_K\}$}
choose $n_1\in\Theta=\{1, \ldots, K\}$ such that $n_1 = \arg\max\limits_{n\in \Theta}\|b_{n}\|$, initialize $\Theta_{in}=\{n_1\}$\; 
compute ${x}_{n_1}$ in \eqref{eq:pPDEs1}, orthogonalize it as basis $w_1$, and initialize $\Phi_r = \left[w_1\right]$\;
\For{$r= 1, \ldots, r_{\max}$}{
1. compute $\widehat{x}_n$ of \eqref{eq:pPDEs2} for all $n \in \Theta\backslash \Theta_{in}$; 
{\footnotesize \tcp{$\CO(m^3K+N^2+N^2K)$ operations at $m$-th iter.}} 
2. evaluate the error estimator $\CE_n$ for all $n \in \Theta\backslash \Theta_{in}$; 
{\footnotesize \tcp{$\CO(mN^2K)$ operations at $m$-th iter.}} 
3. choose $n_{r+1}= \arg\max\limits_{n\in \Theta\backslash \Theta_{in}}\CE_n$\; 
4. \If{$\CE_{n_{r+1}}< {\sf \tol_{ROM}}$} {calculate $x_n = \Phi_r \widehat{x}_n$ for all $n \in \Theta\backslash \Theta_{in}$\; 
break\;}  
5. update $\Theta_{in} = \{\Theta_{in}, n_{r+1}\}$\;
6. compute ${x}_{n_{r+1}}$ to \eqref{eq:pPDEs1} at $n=n_{r+1}$, orthogonalize it as basis $w_{r+1}$, and update $\Phi_r = \left[\Phi_r, w_{r+1}\right]$;  
{\footnotesize\tcp{$\CO(N^2)$ operations at $m$-th iter.}} 
}
\caption{Greedy basis generation}\label{alg:basis}
\end{algorithm}
The major computation inside the \textbf{for} loop lies in steps 1, 2 and 6. At the $m$-th iteration, the complexity of these three steps contains $\CO(m^3K+N^2+N^2K)$ operations for assembling the matrices and solving the ROMs with direct solver; $\CO(mN^2K)$ operations for estimating the residual norm; and $\CO{(N^2)}$ operations for solving the FOM with multigrid V-cycles and obtaining a new basis vector. Therefore, {\em for each FGMRES iteration,} the computation complexity for the Algorithm \ref{alg:basis} with ROM dimension $r\ll N$ is $\CO{(r^4K+2rN^2+rN^2K+r^2N^2K)}=\CO(r^2N^2K) $ operations,
which is theoretically comparable to the multigrid-based FOM solver with the   $\CO(N^2K)$ operations for Step-(b), provided $r$ is bounded by a mesh-independent constant. 
That is, the ROM-based solver requires the underlying dynamics to have a low-dimensional solution manifold, to be computationally competitive with the multigrid-based FOM solver.

In order to further improve the practical computational efficiency of the above greedy basis generation algorithm, various heuristic approaches have been proposed in the past decade, see e.g. in \cite{sen2008reduced,jiang2020adaptive}.
For our greedy basis generation algorithm, we introduce the following algorithmic modifications to reduce its computational complexity.
 

\textbf{(i) Update reduced matrices recursively.} 
The reduced systems can be recursively assembled by making use of the basis structure $\Phi_r=[\Phi_{r-1},w_r]$.
Define $A_n^{(r)}=\Phi_r^* A_n \Phi_r$, $B_r = \Phi_r^*A_n$, $C_r = A_n \Phi_r$, and $F_r = \Phi_r^*S_1$, then we recursively construct them (with $r\ge 2$) according to
\begin{equation*}
A_n^{(r)} = 
\left[
\begin{array}{cc}
A_n^{(r-1)} & B_r w_r \\
w_r^* C_r & w_r^* A_n w_r
\end{array}
\right], \,\, 
B_r =
\left[
\begin{array}{c}
B_{r-1}\\
w_r^*A_n
\end{array}
\right], \,\,
C_r =
\left[
\begin{array}{cc}
C_{r-1} & A_n w_r
\end{array}
\right], \,\, 
F_r =
\left[
\begin{array}{cc}
F_{r-1} & w_r^*S_1
\end{array}
\right],
\end{equation*}
starting with the initialization:
$
	 A_n^{(1)} =  w_1^\intercal A_n w_1 , \,\, 
	B_1 = 
		w_1^*A_n , \,\,
	C_1 =  A_n w_1  , \,\, 
	F_1 = w_1^*S_1 . 
$

\textbf{(ii) Sample $\Theta\backslash \Theta_{in}$ randomly.} 
In steps 1 and 2 at each \textbf{for} loop, instead of testing the entire set $\Theta\backslash \Theta_{in}$ for estimating the ROM approximation errors, we randomly choose $P$ components from it. Define $P=r_p K$ with $r_p\in (0,1]$. Although one usually fixes a randomly sampled training parameter set in RBM, we vary the sample set during the training process in order to achieve generalization. We refer to \cite{cohen2020reduced} for further analysis on RBM using random training sets. Based on our numerical experiments, taking only 10 percent (i.e. $r_p=0.1$) of all the indices in $\Theta\backslash \Theta_{in}$ for training is enough to obtain accurate approximation results. 

\textbf{(iii) Evaluate residual norm on a coarser mesh.} 
In step 2, we compute the residual vectors at $N_s^2$ points instead of all the $N^2$ spatial grid points, with $N_s<N$, in order to reduce the complexity of evaluating the residual norm. 
To define these points, one way is to sample spatial coordinates at random. However, it could become unreliable if the distribution of residuals is far more uneven. In that case, structured random embeddings could be considered \cite{martinsson2020randomized,balabanov2019randomized}. 
Here, we propose to select the points from a coarser spatial mesh. 
Based on our numerical experiments, it is sufficient to use a four times coarser mesh size (i.e. $N_s=\frac{N}{4}$ or  $N_s^2=\frac{N^2}{16}$) for obtaining reliable and accurate results cheaply.

With these modifications (i.e. $P=r_p K$ and $N_s=N/4$), the complexity of Step-(b) {\em at each FGMRES iteration} now becomes $\CO(r^4P+2rN^2+rN^2P+r^2N_s^2P)=
\CO(r_p r^4 K+2rN^2+r_p rN^2K+(r_p/16)r^2 N^2K)$ operations, {which scales as $\CO((r/16+1) r_p r N^2K)$ if $r\ll N$ and hence}
improves the original complexity of Algorithm \ref{alg:basis} since $r_p\in (0,1]$. Therefore, it has a great potential to become more efficient than the multigrid-based FOM solver in practice, especially when $r$ is relatively small.
Such an improvement is conditional but is not surprising because the multigrid algorithm in general delivers the `optimal' complexity for solving elliptic systems.

We also note that the number of outer FGMRES iterations might change when different solvers are used in Step-(b). Based on our numerical tests in next section, the ROM-based solvers in Step-(b) always lead to less number of outer FGMRES iterations than the multigrid-based FOM solver while achieving the same approximation accuracy. 
Finally, we remark that steps 1-2 are parallelizable as individual problems are independent, the involved matrix-vector products can also be implemented in a parallel fashion. One can further divide all the $K$ systems into disjoint subsets and then run the Algorithm 1 within each subset. Some discussions on parallel computing in the RBM setting can be found in \cite{knezevic2011high}.

\subsection{The motivating example revisited}

\begin{figure}[!htb]
	\begin{center}
		\includegraphics[width=0.48\textwidth]{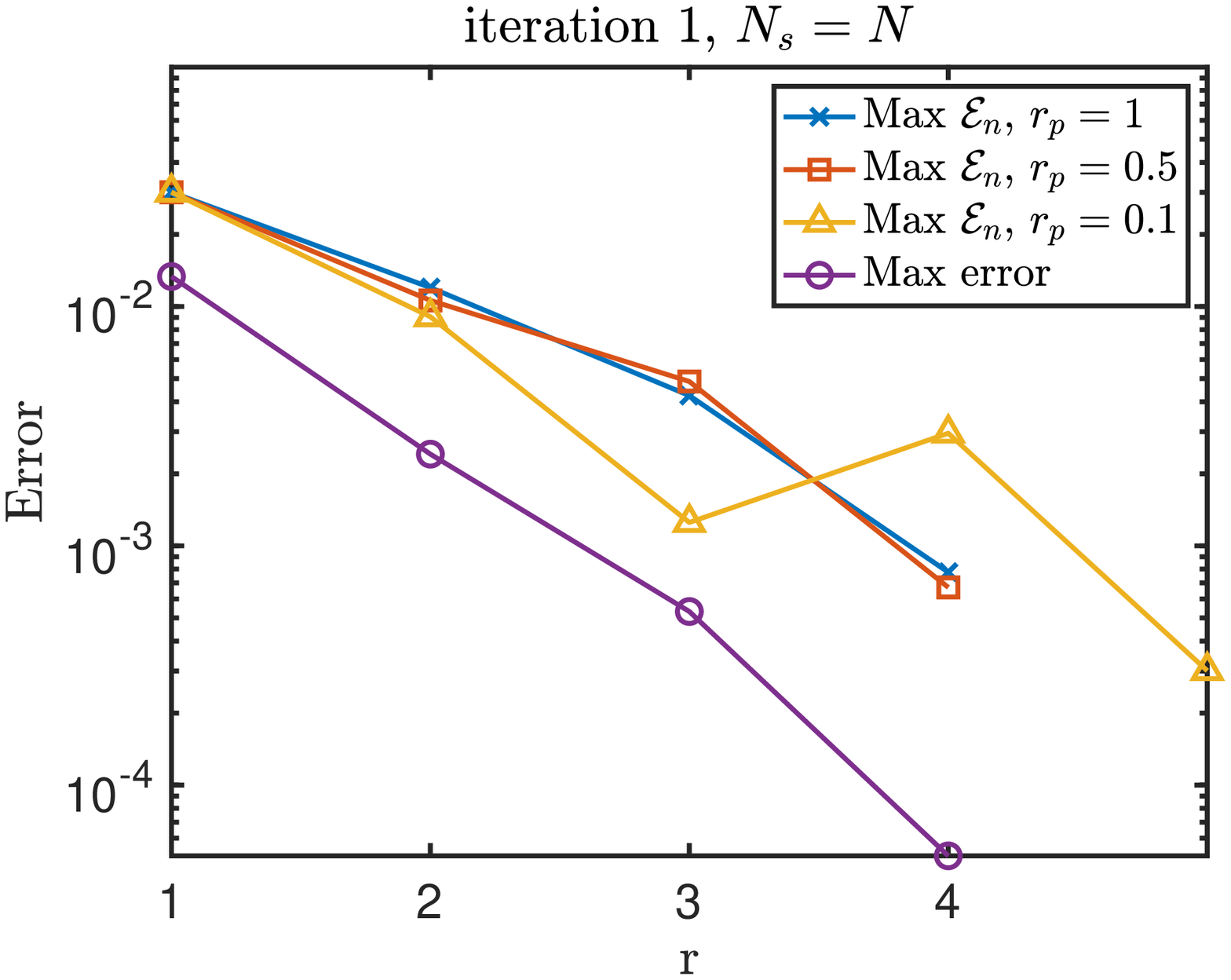} \quad
		\includegraphics[width=0.48\textwidth]{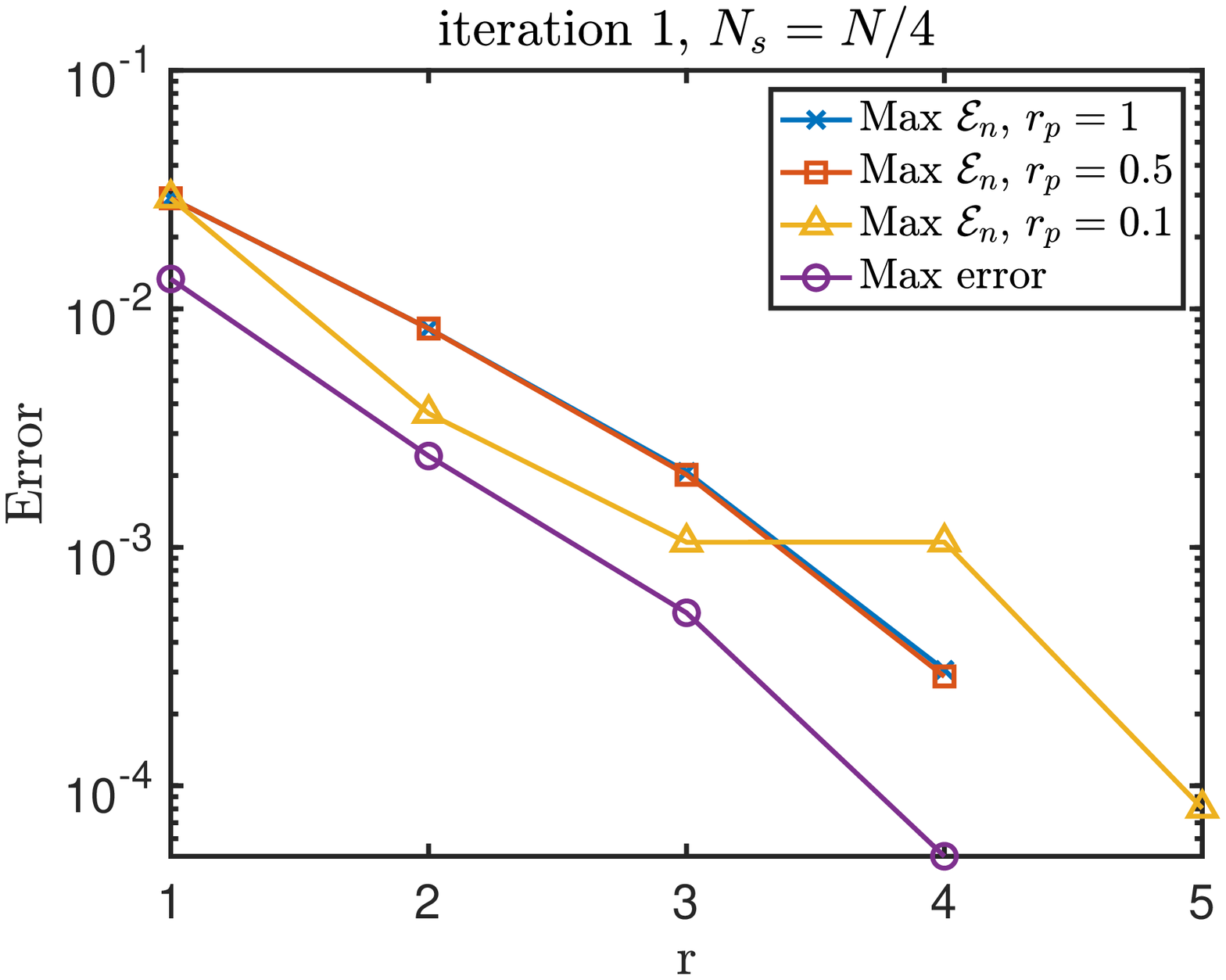} \\
		\includegraphics[width=0.48\textwidth]{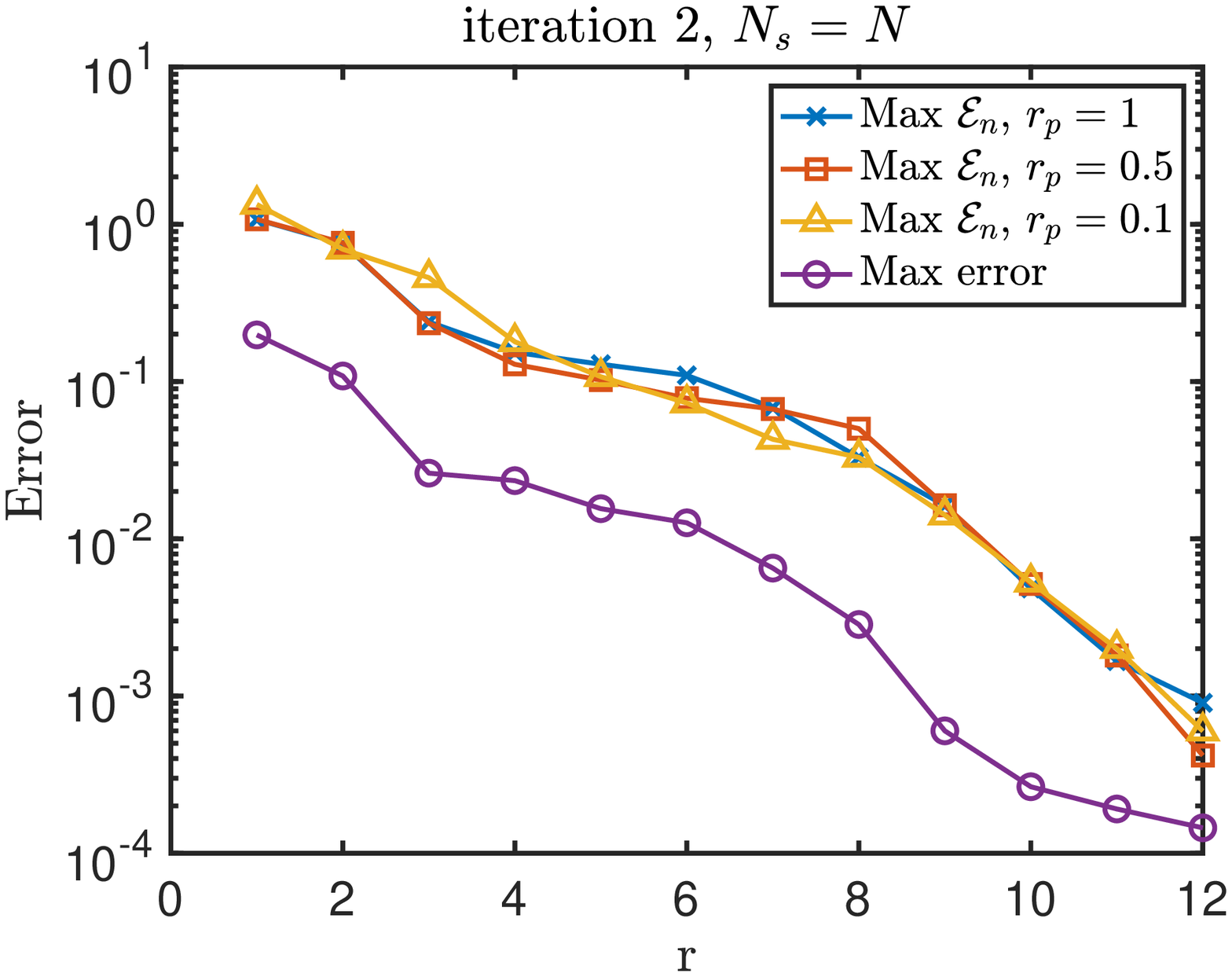} \quad
		\includegraphics[width=0.48\textwidth]{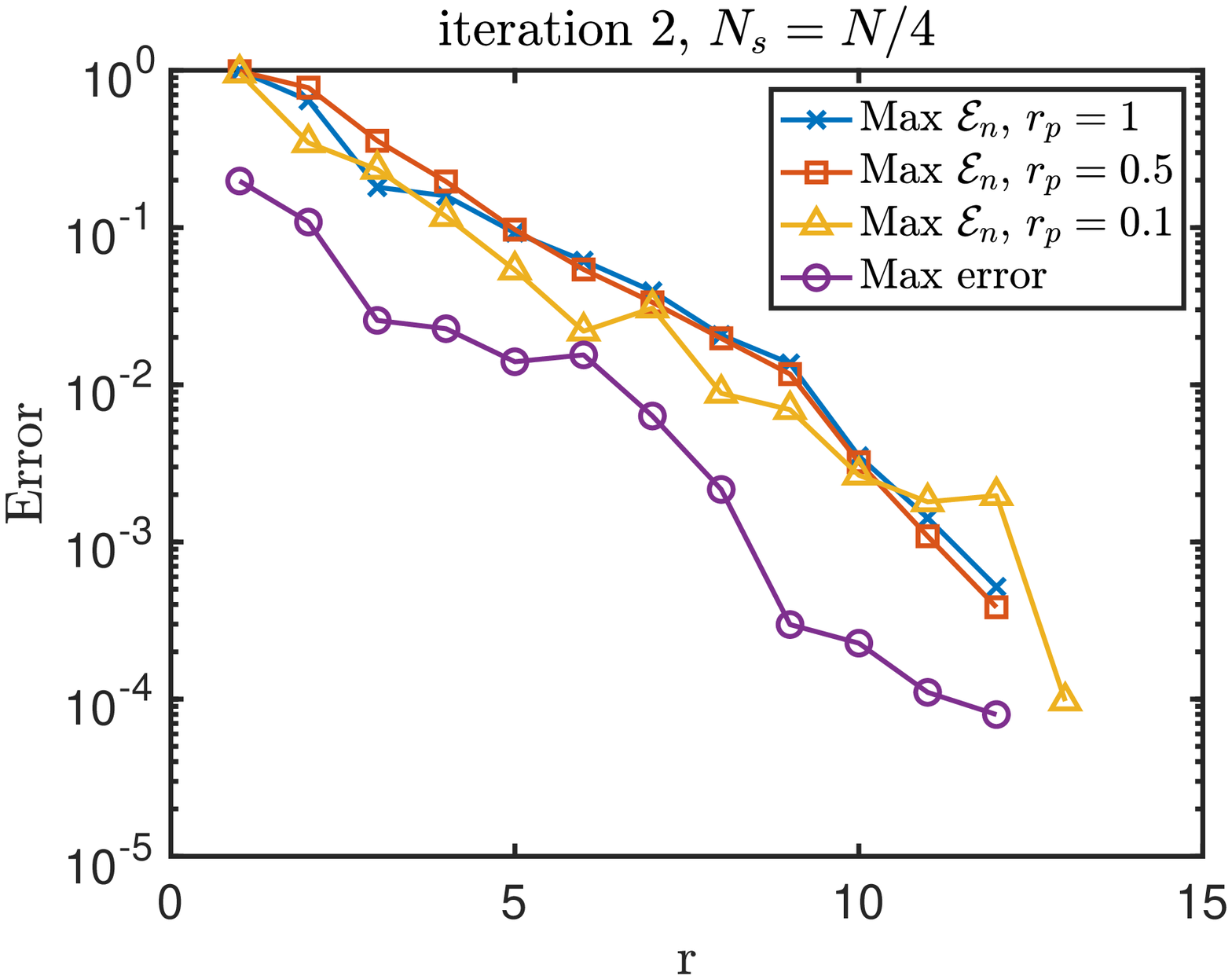}
	\end{center}
	\caption{Exact errors and estimated errors of Step-(b) at iteration 1 (first row) and iteration 2 (second row) while varying the value of $r_p= 1, 0.5, 0.1$: (left) $N_s=N$ and (right) $N_s=\frac{N}{4}$.  
	}
	\label{plot_2D_ex1_err}
\end{figure}
We illustrate the error indicator $\CE_n$ of our proposed ROM-based solver in Step-(b) by investigating the same motivating example again. 
The aforementioned modifications for Algorithm \ref{alg:basis} involve the use of random training sets and evaluations of residual norms on a coarser spatial mesh. 
First, under the same discretization as in Section \ref{sec:me1}, we keep $N_s=N$ but vary $r_p$ from $1$, $0.5$ to $0.1$. In these cases, the preconditioned FGMRES converges in 2 iterations. 

Fig. \ref{plot_2D_ex1_err} (left column) shows the evolution of estimated errors with respect to the increasing number of reduced basis vectors during the greedy search algorithm, together with the exact errors. It is seen that: (1) the error indicator gives a good approximation of the actual approximation error; (2) with the same stopping criterion, the dimension of resulting reduced basis stays almost the same when different values of $r_p$ (or equivalently $P$) are used. Therefore, as $r_p$ gets smaller, the computational efficiency improves but the overall accuracy of the algorithm does not deteriorate. Using the same computational setting but decreasing $N_s$ from $N$ to ${N}/{4}$, we observe similar numerical behaviors as shown in Fig. \ref{plot_2D_ex1_err} (right column). This indicates the used heuristic modifications for Algorithm \ref{alg:basis} are practicable and effective. 

 \begin{figure}[!htb]
 	\begin{center}
 		\includegraphics[width=0.48\textwidth]{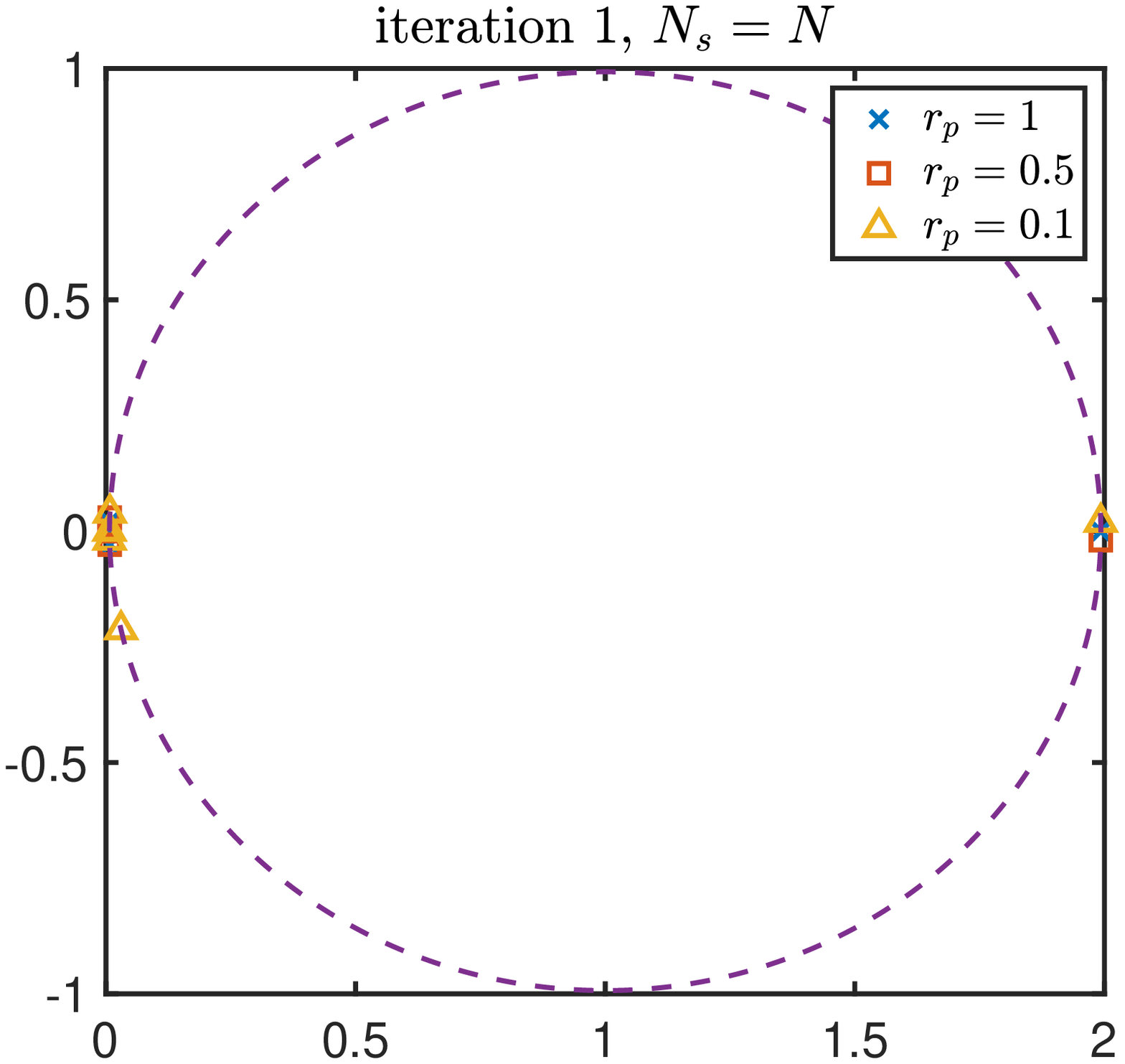} \quad
 		\includegraphics[width=0.48\textwidth]{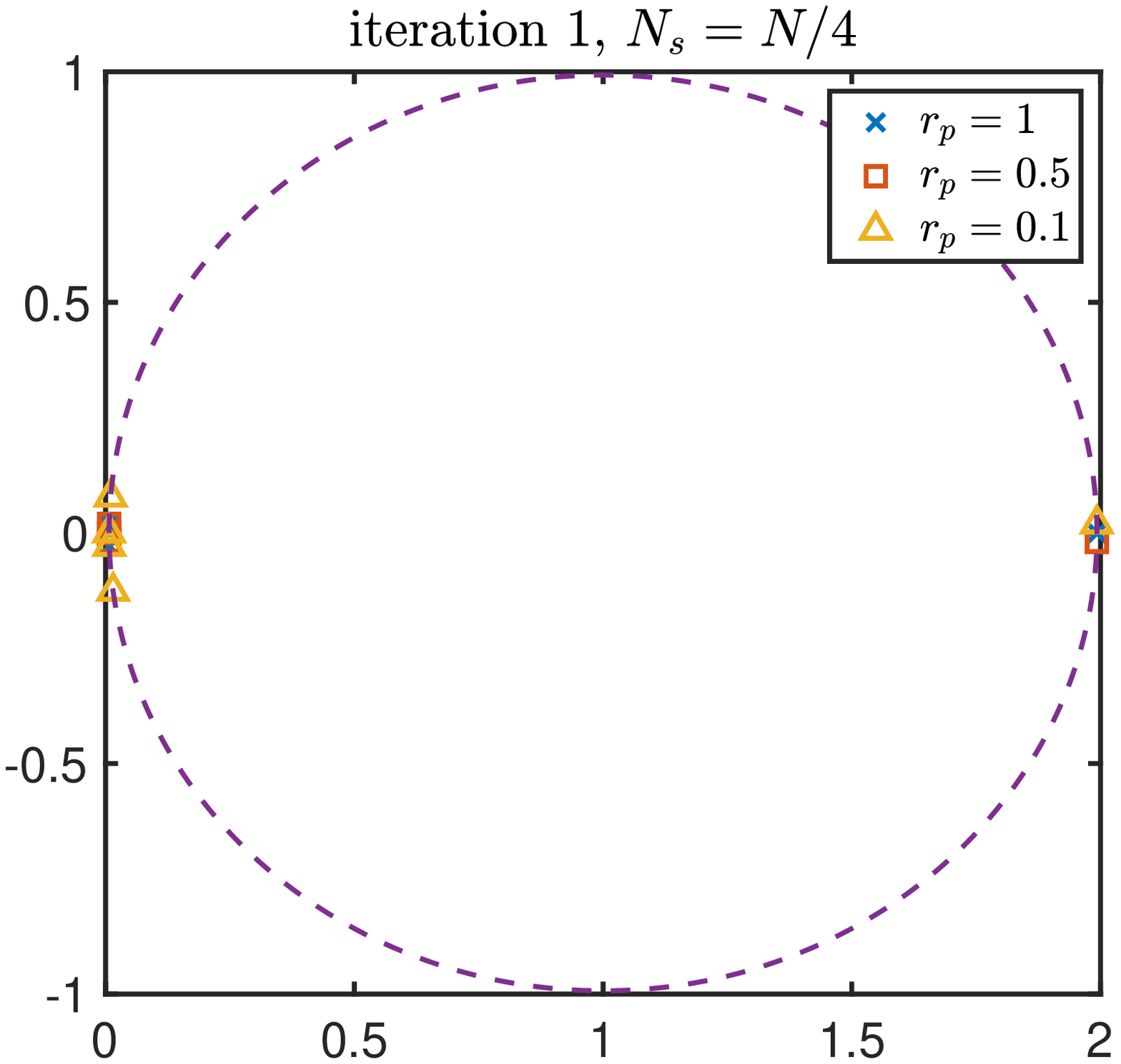} \\
 		\includegraphics[width=0.48\textwidth]{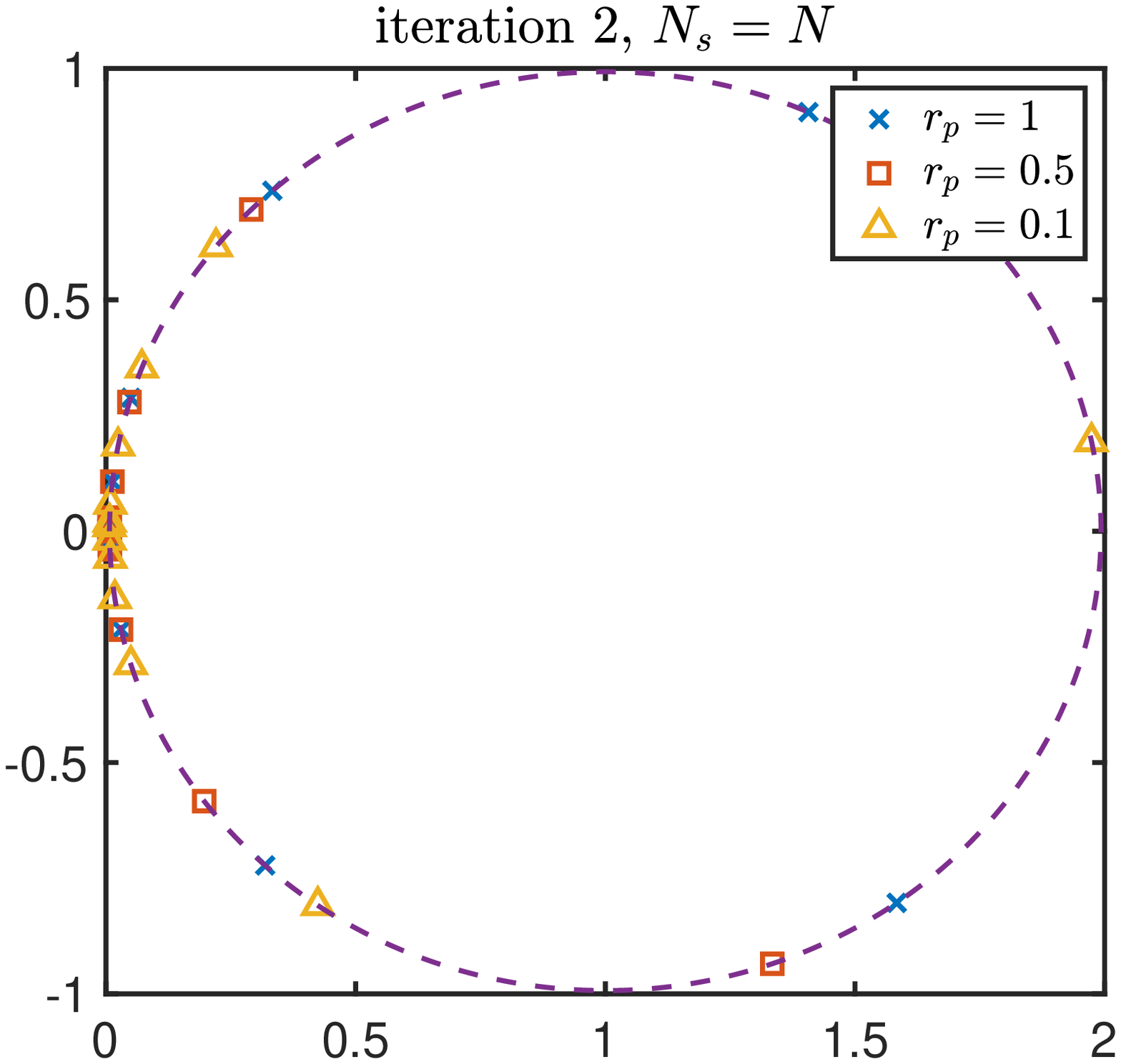} \quad
 		\includegraphics[width=0.48\textwidth]{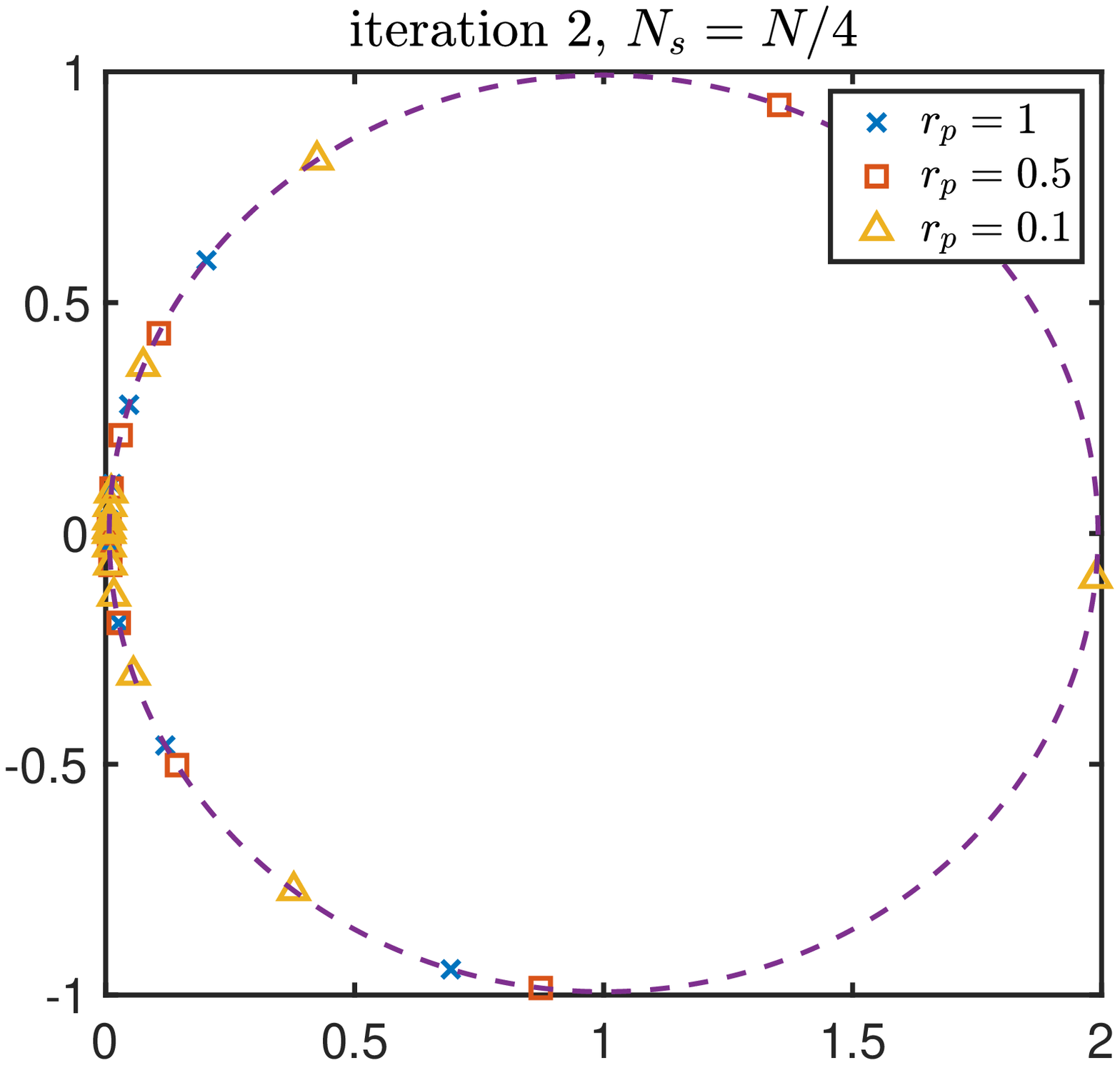}
 	\end{center}
 	\caption{Distribution of $d_n$ for $n\in\Theta_{in}$ in Step-(b) at iteration 1 (first row) and iteration 2 (second row) while varying the value of $r_p= 1, 0.5, 0.1$: (left) $N_s=N$ and (right) $N_s=\frac{N}{4}$.  
 	}
 	\label{plot_2D_ex1_dn}
 \end{figure}
We plot the complex number $d_n$'s associated to the selected basis indices $n\in \Theta_{in}$ as points in Fig. \ref{plot_2D_ex1_dn}. The points distribute on the  circle centered at $1+0\i$ with radius $\alpha^{\frac{1}{K}}\approx 0.99$ but concentrate near $\theta_n=0$ or $2\pi$. 
It matches our observations in Fig. \ref{plot_2D_ex1_S2} where more structures appear near $n=1$ and $n=K$. 
This also implies ROMs generated from evenly distributed snapshots (as commonly used in POD-based ROMs) may not be effective enough for such systems.  

Compared with the MG-based FOM solver with 4 FGMRES iterations, the proposed ROM-based solver, when $r_p=0.1$ and $N_s=\frac{N}{4}$, achieves the same accuracy $5.7\times 10^{-5}$ with only 2 FGMRES iterations. As expected, the total CPU time decreases about 9 times from 6.2~s to 0.7~s, in which 0.4~s (about 57\%) is spent on solving Step-(b). The reduced number of preconditioned FGMRES iterations also contributes to the substantial saving of overall CPU times.

	\section{Numerical examples}\label{secNum}
	In this section, we present several 1D and 2D PDE examples to illustrate the effectiveness of our proposed ROM-based preconditioning techniques.
	All simulations are implemented using MATLAB on a Dell Precision 7520 Laptop
	with Intel(R) Core(TM) i7-7700HQ CPU@2.80GHz and 48GB RAM.
	The CPU time (in seconds) for solving the whole system and all the sub-systems in Step-(b) are estimated separately by using the timing functions \texttt{tic}/\texttt{toc} with our serial implementation.
	We employ the right-preconditioned FGMRES \cite{saad1993flexible} solver (without restarts) provided in the Tensor Train (TT) Toolbox \footnote{\url{https://github.com/oseledets/TT-Toolbox}} 
	and choose a zero initial guess with a stopping tolerance ${\tol}=10^{-6}$ based on the reduction in relative residual norms.
	To avoid introducing possible round-off errors due to $\alpha$ being too small (as theoretically preferred), we will choose the preconditioner $P_\alpha$ with a fixed small $\alpha=10^{-2}$, which leads to only a small number of preconditioned FGMRES iterations in all tested examples.
	
 	For the ROM-based solver, we choose the relative residual tolerance to be ${\sf\tol_{ROM}=\sqrt{\tol}}$, which works well for all the tested examples. 	In particular, we numerically observed that using a costly sparse direct solver in Step-(b) yields the same outer FGMRES iteration numbers as our ROM-based preconditioner (ParaDIAG-ROM). We highlight that a smaller ROM residual tolerance would lead to higher reduced basis dimension and computational costs, which however will not improve the overall FGMRES accuracy. 
Within the Algorithm 1, we will solve the selected full-order systems \eqref{eq:pPDEs1} by the backslash sparse direct solver for 1D examples (with tridiagonal systems), and the geometric multigrid V-cycles with ILU smoother and the same stopping tolerance ${\tol}=10^{-6}$ for 2D examples.
 	We set the random number generator using the function \textsf{rng(0,'v5uniform')}, with a seed of 0 and the uniform generator \textsf{'v5uniform'}.

 For the compared multigrid-based ParaDIAG (ParaDIAG-MG) preconditioner, only one V-cycle is used to approximately solve inner linear systems in Step-(b), which gives somewhat rough approximation and therefore requires more outer FGMRES iterations than the ROM-based preconditioner with our chosen tolerance. 
 In ParaDIAG-MG preconditioner, for a general comparison purpose we will use the Gauss-Seidel (GS) smoother and incomplete LU(ILU) smoother for 1D and 2D examples, respectively. For 2D examples, the ILU smoother usually provides more robust and faster convergence rate than the GS smoother. 
 We have tried to vectorize the majority of our MATLAB codes whenever possible.

In all the numerical experiments, we compare the errors and CPU time between the ParaDIAG-MG and ParaDIAG-ROM. 
	The approximation errors are evaluated in discrete $L^2$ norm in both space and time. The order of accuracy is estimated by the logarithmic ratio of the approximation errors between two successive refined (halved) meshes, which should be close to 1 for the numerical scheme if solutions are sufficiently smooth. 
	The CPU time measures the efficiency of the two approaches in serial computing. Since the only difference between them is the solver used in Step-(b), we also list the CPU time for this particular step (displayed inside parentheses). 
	
  \subsection{1D Examples.}
  
\begin{figure}[!htb]
	\begin{center}
	\includegraphics[width=1\textwidth]{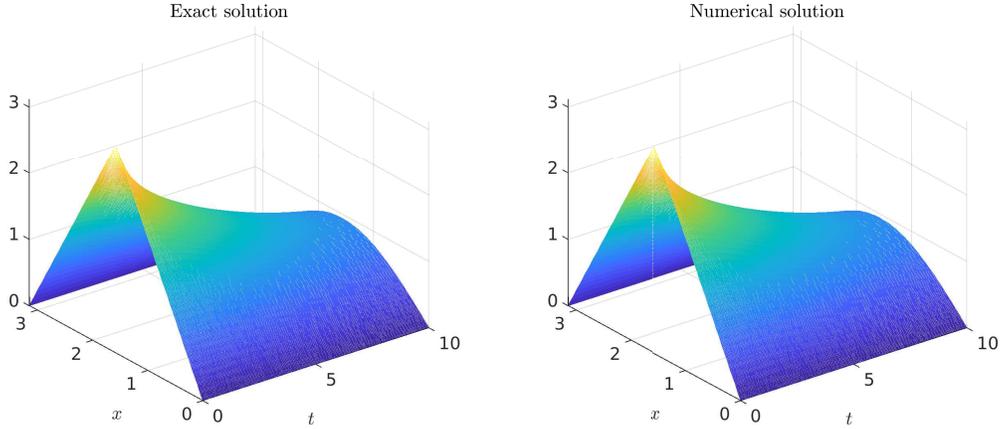} 
	\end{center}
	\caption{The exact solution and computed numerical solution by ParaDIAG-ROM preconditioner in Example 1a.
	}
	\label{plot_1D_ex1}
\end{figure}

\textbf{Example 1a. Heat equation with constant diffusion coefficients.}   We first consider
  a linear heat equation with constant coefficients {defined on the domain}  $\Omega=(0,\pi)$:
  \begin{equation}\label{heatPDE1D}
  	\begin{cases}
  		u_t= \epsilon u_{xx} +f, &\tn{in} \Omega\times(0,T),\\
  		u=0,&\tn{on} \partial\Omega\times(0, T).
  	\end{cases}
  \end{equation} 
  We choose $f=0$, and a non-smooth triangle shaped initial condition 
  \[u(x,0)=\begin{cases}
  	2x, &  0\le x\le \pi/2, \\
  	2(\pi-x), &  \pi/2<x \le \pi,
  \end{cases}
  \]
  with the exact solution is given by
  \[
  u(x,t)=\frac{8}{\pi}\sum_{n=0}^{\infty}  \frac{\cos((2n+1) (2 x-\pi)/2)}{(2n+1)^2} e^{-\epsilon (2n+1)^2 t}.
  \]
  The errors and convergence results are reported in Table \ref{Tab_1D_ex1},
  where the total CPU time of our proposed ROM-based ParaDIAG (ParaDIAG-ROM) preconditioner
  is over 10 times faster than the multigrid-based ParaDIAG (ParaDIAG-MG) preconditioner with the point-wise Gauss-Seidel smoother. 
  The speed up ratio in CPU time for computing Step-(b) is even higher than 20 times.
  Due to the non-smoothness of the initial condition, the order of accuracy is expected to be slightly lower than one as the mesh refines.
  Since one multigrid V-cycle very approximately solve the linear systems in Step-(b), which leads to more outer FGMRES iterations than our ParaDIAG-ROM preconditioner.  
  The average  dimension of the ROM reduced basis, denoted by column `$\bar r$', shows very mild growth as the mesh is refined, which is reasonable considering the full order models' dimension increases by four times.
  Fig. \ref{plot_1D_ex1} compares the exact solution and the numerical solution computed by  our ParaDIAG-ROM preconditioner,
  where the non-smooth initial condition is smoothed out quickly due to the diffusion.
   	\begin{table}[!htb]
  	\centering 
  	\caption{Results of preconditioned FGMRES  for Example 1a: 1D heat equation ($\epsilon=0.1$, $T=10$)}
  	\begin{tabular}{|c||cccc||ccccc|ccc|cccc|}\hline
  		&\multicolumn{4}{c||}{ParaDIAG-MG(GS) Preconditioner}  &\multicolumn{5}{c|}{ParaDIAG-ROM Preconditioner}
  		\\
  		\hline
  		$(N,K)$& Error & Order  & Iter &CPU& Error & Order  & Iter  &CPU & $\bar r$\\   \hline 
(256,2560)&	8.5E-04&	1.2 &	 6&	 4.1 (4.0)&       8.5E-04&	1.2 &	 2&	 0.2 (0.1)&	10\\
(512,5120)&	4.1E-04&	1.1 &	 6&	 12.5 (11.6)&	  4.1E-04&	1.1 &	 2&	 0.9 (0.5)&	12\\
(1024,10240)&	2.2E-04&	0.9 &	 6&	 42.3 (38.1)	&	2.2E-04&	0.9 &	 2&	 3.3 (1.8)&	13\\ 
  		\hline
  	\end{tabular}
  	\label{Tab_1D_ex1}
  \end{table}

\begin{figure}[!htb]
  	\begin{center}
  		\includegraphics[width=1\textwidth]{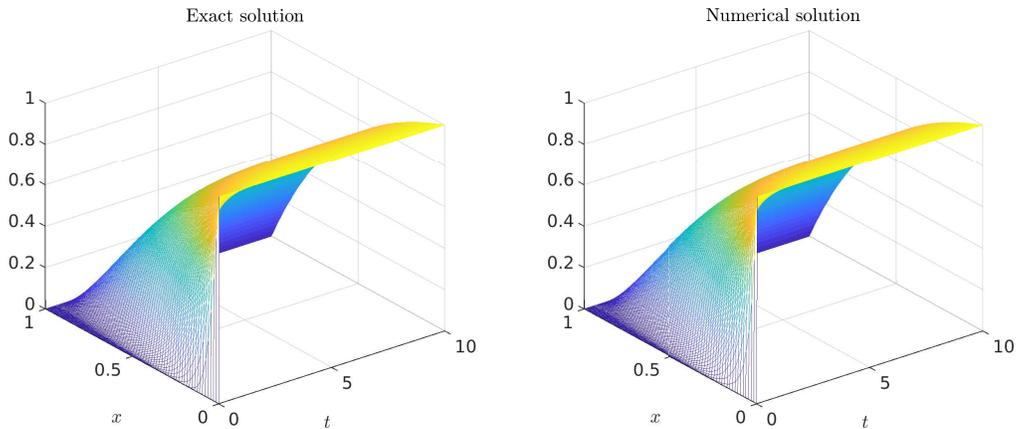} 
  	\end{center}
  	\caption{The exact solution and computed numerical solution by ParaDIAG-ROM preconditioner in Example 1b.
  	}
  	\label{plot_1D_ex4}
  \end{figure}

  \textbf{Example 1b. Convection-diffusion (C-D) equation.}  
  We also consider
  another convection-diffusion equation  {defined on the domain}  $\Omega=(0,1)$:
  \begin{equation}\label{CD1D_a}
  	\begin{cases}
  		u_t= \epsilon u_{xx}-c u_x, &\tn{in} \Omega\times(0,T),\\
  		u(0,t)=1,\quad u(1,t)=0,&  t\in(0, T),\\
  		u(x,0)=0,&  x\in \Omega,
  	\end{cases}
  \end{equation} 
  where $P_e=c/\epsilon$ is the P\'{e}clet number that represents the ratio of the convection rate over the diffusion rate. The exact solution is explicitly known \cite{mohsen1983analytical} in terms of infinity series, which is used to measuring the errors. 
  The errors and convergence results are reported in Table \ref{Tab_1D_ex4},
  where our ParaDIAG-ROM preconditioner
  is over 10 times faster than the ParaDIAG-MG preconditioner. Fig. \ref{plot_1D_ex4} compares the exact solution and the numerical solution computed by  our ParaDIAG-ROM preconditioner. 
  Due to a discontinuity at the corner $(0,0)$, our used finite difference scheme show a slightly degraded order of accuracy.
  Again, the boundary layer leads to slightly larger  $\bar r$ for accurate approximation than the previous example.
  \begin{table}[!htb]
  	\centering 
  	\caption{Results of preconditioned FGMRES  for Example 1b: 1D C-D equation with upwind scheme and ($\epsilon=0.1, c=0.2$, $T=10$)}
  	\begin{tabular}{|c||cccc||ccccc|ccc|cccc|}\hline
  		&\multicolumn{4}{c||}{ParaDIAG-MG(GS) Preconditioner}  &\multicolumn{5}{c|}{ParaDIAG-ROM Preconditioner}
  		\\
  		\hline
  		$(N,K)$& Error & Order  & Iter &CPU& Error & Order  & Iter  &CPU & $\bar r$\\   \hline 
(256,2560)&	6.1E-03&	0.7 &	 5&	 3.5 (3.4)&	       6.1E-03&	0.7 &	 2&	 0.2 (0.1)&	13\\
(512,5120)&	3.7E-03&	0.7 &	 5&	 10.1 (9.4)&	   3.7E-03&	0.7 &	 2&	 0.7 (0.4)&	14\\
(1024,10240)&	2.3E-03&	0.7 &	 5&	 35.5 (32.1)	&	2.2E-03&	0.7 &	 2&	 2.9 (1.5)&	16\\
  		\hline
  	\end{tabular}
  	\label{Tab_1D_ex4}
  \end{table}
  
  \subsection{2D Examples.}
\begin{figure}[!htb]
 	\begin{center}
 		\includegraphics[width=1\textwidth]{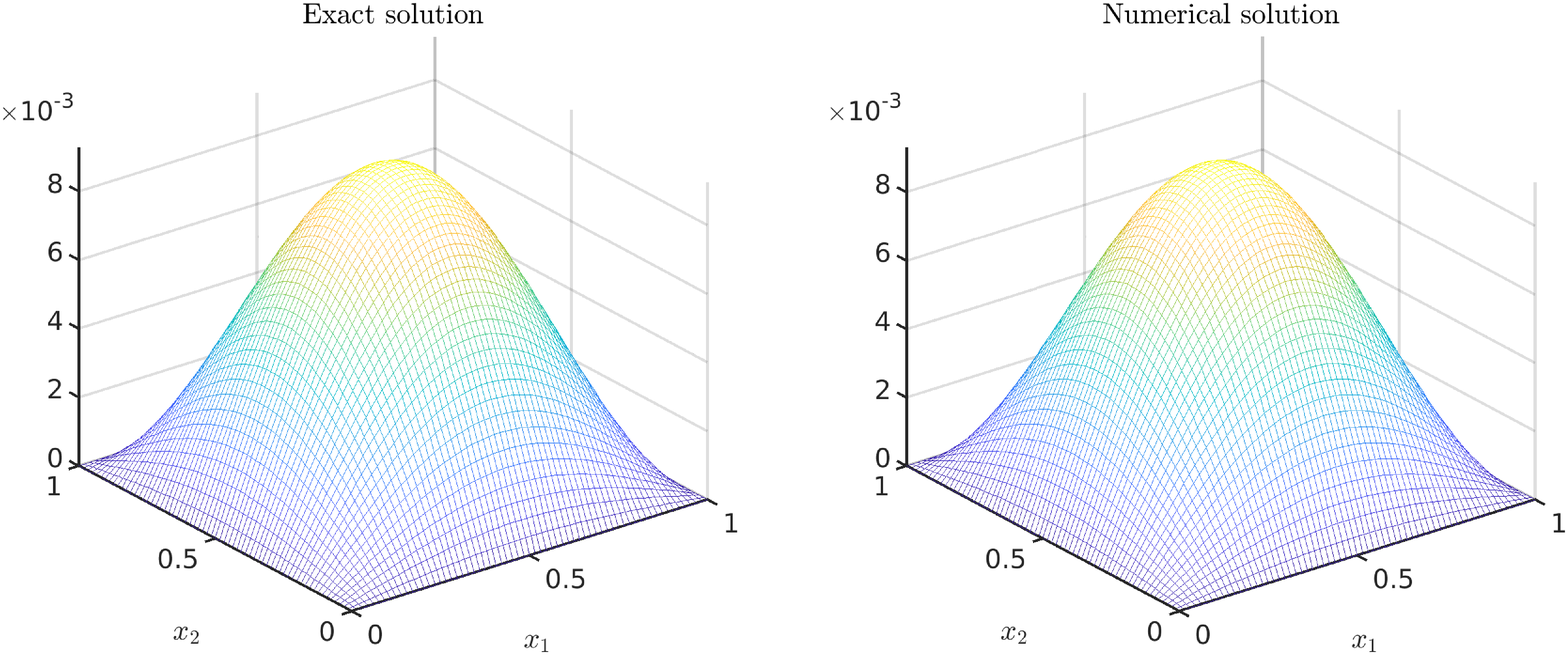} 
 	\end{center}
 	\caption{The exact solution and computed numerical solution by ParaDIAG-ROM preconditioner at the final time $T=10$ in Example 2a.
 	}
 	\label{plot_2D_ex1}
 \end{figure}

  \textbf{Example 2a. Heat equation with constant diffusion coefficients.}
We consider the motivating example used in Section \ref{sec:me1} again. 
 The errors and convergence results are reported in Table \ref{Tab_2D_ex1},
where our ParaDIAG-ROM preconditioner
is about 8 times faster than the ParaDIAG-MG preconditioner. 
In the last column `$\bar r$ (Vc)', the first number $\bar r$ indicates the average reduced basis dimension and Vc stands for the average multigrid V-cycles used for solving all the chosen full-order models.
The ParaDIAG-MG preconditioner shows very fast and robust convergence rate due to the used ILU smoother, but the overall computational cost is still high since it requires total $K$ V-cycles to approximately solve all the $K$ linear systems at each iteration.
In contrast, the ParaDIAG-ROM preconditioner only used in average $\bar r\times \tn{(Vc)}\ll K$ full-order model V-cycles to find the reduced basis and then solved $K$ reduced linear systems of much smaller sizes cheaply.
Fig. \ref{plot_2D_ex1} compares the exact solution and the numerical solution  computed by  our ParaDIAG-ROM preconditioner at the final time. 

	 \begin{table}[!htb]
	 	\centering 
	 	\caption{Results of preconditioned FGMRES for Example 2a: 2D heat equation ($\epsilon=0.01$, $T=10$)}
	 	\begin{tabular}{|c||cccc||ccccc|ccc|cccc|}\hline
	 		&\multicolumn{4}{c||}{ParaDIAG-MG(ILU) Preconditioner}  &\multicolumn{5}{c|}{ParaDIAG-ROM Preconditioner}
	 		\\
	 		\hline
	 		$(N^2,K)$& Error & Order  & Iter &CPU& Error & Order  & Iter  &CPU & $\bar r$ (Vc)\\   \hline 
($64^2$,640)&		 5.7E-05&	 1.2 &	 4&	 6.2(5.5)&	      5.7E-05&	 1.2 &	 2&	 0.7 (0.4)&	8 (4) \\ 
($128^2$,1280)&		 2.7E-05&	 1.1 &	 4&	 48.3(40.7)&	  2.7E-05&	 1.1 &	 2&	 6.1 (2.6)&	9 (5) \\ 
($256^2$,2560)&		 1.3E-05&	 1.0 &	 4&	 378.5(301.5)&	  1.3E-05&	 1.0 &	 2&	 43.1 (20.3)&	10 (5) \\ 
	 		\hline
	 	\end{tabular}
	 	\label{Tab_2D_ex1}
	 \end{table}
 
\begin{figure}[!htb]
  	\begin{center}
  		\includegraphics[width=1\textwidth]{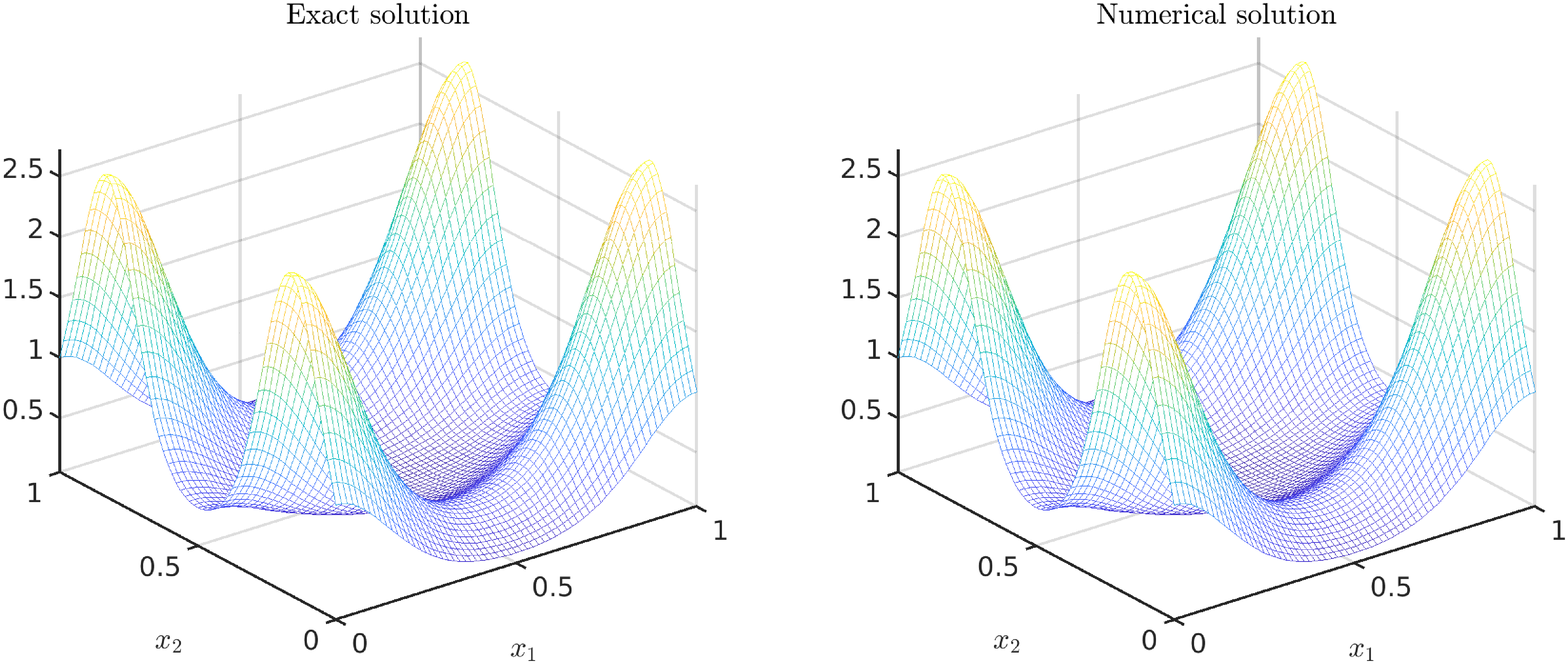} 
  	\end{center}
  	\caption{The exact solution and computed numerical solution by ParaDIAG-ROM preconditioner at the final time $T=10$ in Example 2b.
  	}
  	\label{plot_2D_ex2}
  \end{figure}

 \textbf{Example 2b. Heat equation with spatially variable diffusion coefficients.}
 We also consider linear heat equation with variable diffusion coefficients (in divergence form) on $\Omega=(0,1)^2$:
 \begin{equation}\label{heatPDEv}
 	\begin{cases}
 		u_t= \nabla\cdot(a\nabla u) +f, &\tn{in} \Omega\times(0,T),\\
 		u=g,&\tn{on} \partial\Omega\times(0, T). 
 	\end{cases}
 \end{equation} 
 We choose $a(x,y)=10^{-5} \sin(\pi x y)$, the Dirichlet boundary condition $g$, the initial condition $u_0$
   and a suitable $f$ such that the exact solution is given as
$ u(x,y,t)=e^{-t/10} e^{\cos(2\pi x)+\sin(3\pi y)}.
$
 The errors and convergence results are reported in Table \ref{Tab_2D_ex2},
 where our ParaDIAG-ROM preconditioner
 is about 10 times faster than the ParaDIAG-MG preconditioner. 
 In particular, it takes only two reduced basis vectors to accurately approximate the full-order model solutions in Step-(b). 
 The variable diffusion coefficients does not affect the convergence rate of both preconditioners.
Our ParaDIAG-ROM preconditioner show very high computational efficiency,
 although the ParaDIAG-MG preconditioner based on ILU smoother indeed converges faster than the 1D case based on the GS smoother. Fig. \ref{plot_2D_ex2} compares the exact solution and the numerical solution  computed by  our ParaDIAG-ROM preconditioner at the final time. 
 \begin{table}[!htb]
 	\centering 
 	\caption{Results of preconditioned FGMRES  for Example 2b: 2D heat equation with variable coefficient ($T=10$)}
 	\begin{tabular}{|c||cccc||ccccc|ccc|cccc|}\hline
 		&\multicolumn{4}{c||}{ParaDIAG-MG(ILU) Preconditioner}  &\multicolumn{5}{c|}{ParaDIAG-ROM Preconditioner}
 		\\
 		\hline
 		$(N^2,K)$& Error & Order  & Iter &CPU& Error & Order  & Iter  &CPU & $\bar r$(Vc)\\   \hline 
($64^2$,640)&		 2.5E-03&	 1.0 &	 3&	 4.4(3.8)&	        2.5E-03&	 1.0 &	 2&	 0.5 (0.1)&	2 (2) \\ 
($128^2$,1280)&		 1.3E-03&	 1.0 &	 3&	 35.9(29.7)&        1.3E-03&	 1.0 &	 2&	 3.7 (0.8)&	2 (2) \\ 
($256^2$,2560)&		 6.5E-04&	 1.0 &	 3&	 289.7(223.2)&	    6.3E-04&	 1.0 &	 2&	 27.0 (6.4)&	2 (2) \\ 
 		\hline
 	\end{tabular}
 	\label{Tab_2D_ex2}
 \end{table}

\begin{figure}[!htb]
	\begin{center}
		\includegraphics[width=1\textwidth]{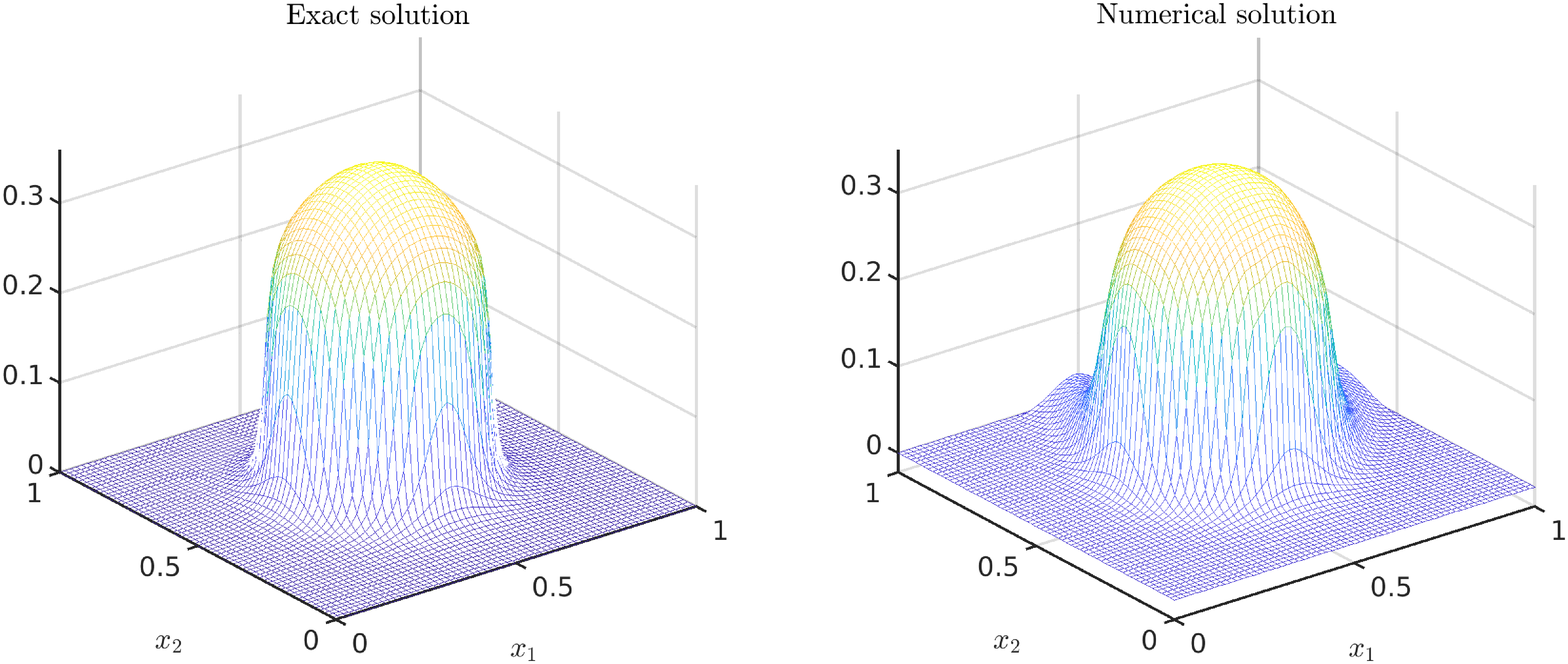} 
	\end{center}
	\caption{The exact solution and computed numerical solution by ParaDIAG-ROM preconditioner at the final time $T=10$ in Example 2c.
	}
	\label{plot_2D_ex3}
\end{figure}

 \textbf{Example 2c. Convection-diffusion (C-D) equation with internal layer.}
  We consider
 another convection-diffusion equation (adapted from \cite{john2008finite}) {defined on the domain}  $\Omega=(0,1)^2$:
 \begin{equation}\label{CD2Db}
 	\begin{cases}
 		u_t=\epsilon \Delta u -\bm c\cdot  \nabla u+f  , &\tn{in} \Omega\times(0,T),\\
 		u(x,y,t)=0,&\tn{on} \partial\Omega\times(0, T)\\
 		u(x,y,0)=u_0(x,y), & \tn{in} \Omega,
 	\end{cases}
 \end{equation} 
where $\bm c=(2,3)$, $u_0=u(x,y,0)$, and $f$ is chosen such that the exact solution reads
\[
u=16e^{-t}x(1-x)y(1-y)\left(1/2+
\arctan(2\sqrt{1/\epsilon}(0.25^2-(x-0.5)^2-(y-0.5)^2))/\pi\right).
\]
 In this problem, there is a hump changing its height over time and the steepness (or thickness) of the circular internal layer depends on the diffusion coefficient $\epsilon$ (or $\sqrt{\epsilon}$).
 The errors and convergence results are reported in Table \ref{Tab_2D_ex3},
 where our ParaDIAG-ROM preconditioner
 is about 3 times faster than the ParaDIAG-MG preconditioner. Fig. \ref{plot_2D_ex3} compares the exact solution and the numerical solution  computed by  our ParaDIAG-ROM preconditioner at the final time. 
 As also observed in \cite{john2008finite}, 
 our used finite difference upwind scheme indeed leads to spurious oscillations behind the hump in the direction of the convection, which
 may be further suppressed by other specially designed schemes.
In this example, the ParaDIAG-MG preconditioner based on ILU smoother converges very fast in only 2-3 iterations, while our ParaDIAG-ROM preconditioner shows slightly increasing $\bar r$(V-cycles), which indicates the advantage of our ParaDIAG-ROM preconditioner may become slightly less significant.
This drawback of mildly increasing average reduced basis   dimension $\bar r$ in our ParaDIAG-ROM preconditioner for convection-dominated problems deserves further investigation.
 \begin{table}[!htb]
 	\centering 
 	\caption{Results of preconditioned FGMRES  for Example 2c: 2D C-D equation with upwind scheme ($\epsilon=10^{-4}$, $T=10$)}
 	\begin{tabular}{|c||cccc||ccccc|ccc|cccc|}\hline
 		&\multicolumn{4}{c||}{ParaDIAG-MG(ILU) Preconditioner}  &\multicolumn{5}{c|}{ParaDIAG-ROM Preconditioner}
 		\\
 		\hline
 		$(N,N,K)$& Error & Order  & Iter &CPU& Error & Order  & Iter  &CPU & $\bar r$ (Vc)\\   \hline 
($64^2$,640)&		 1.0E-01&	 0.8 &	 2&	 3.0(2.6)&	       1.0E-01&	 0.8 &	 2&	 0.8 (0.4)&	13 (2) \\ 
($128^2$,1280)&		 6.1E-02&	 0.7 &	 2&	 24.8(19.9)&	   6.1E-02&	 0.7 &	 2&	 5.8 (2.8)&	17 (2) \\ 
($256^2$,2560)&		 3.6E-02&	 0.8 &	 3&	 285.9(221.9)&	   3.6E-02&	 0.8 &	 2&	 46.9 (24.8)&	21 (3) \\ 
 		\hline
 	\end{tabular}
 	\label{Tab_2D_ex3}
 \end{table}

\begin{figure}[!htb]
 	\begin{center}
 		\includegraphics[width=1\textwidth]{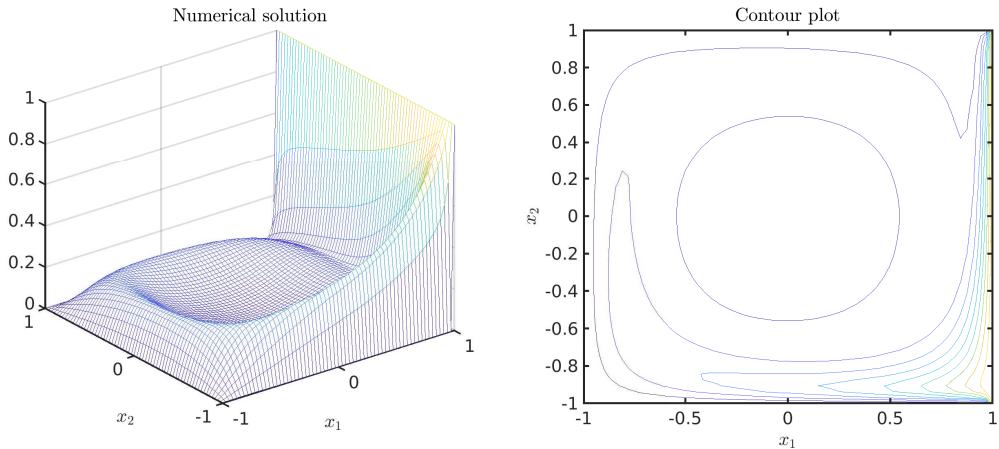} 
 	\end{center}
 	\caption{The computed numerical solution by ParaDIAG-ROM preconditioner and its contour plot with $N=64$  at the final time $T=10$ in Example 2d.
 	}
 	\label{plot_2D_ex4}
 \end{figure}

\textbf{Example 2d. Convection-diffusion (C-D) equation with boundary layer.}
In the last example, we consider
the prototype convection-diffusion equation \cite{ESW05,lin2020all} {defined on the domain}  $\Omega=(0,1)^2$:
\begin{equation}\label{CD2D}
	\begin{cases}
		u_t=\epsilon \Delta u -\bm c\cdot  \nabla u, &\tn{in} \Omega\times(0,T),\\
		u(x,y,t)=(1-e^{-10t})\mathbbm{1}_{\{x=1\}},&\tn{on} \partial\Omega\times(0, T)\\
		u(x,y,0)=0 &\tn{in} \Omega.
	\end{cases}
\end{equation} 
where $\bm c=(2y(1-x^2),-2x(1-y^2))$  is the circulating wind
and the indicator function $\mathbbm{1}_{\{x=1\}}$ denotes one-side hot wall boundary condition at $x=1$. 
Since the  exact solution is unknown, we will only report the relative residual norms of the computed solutions.
The errors and convergence results are reported in Table \ref{Tab_2D_ex4},
where our ParaDIAG-ROM preconditioner
is about 8 times faster than the ParaDIAG-MG preconditioner. Fig. \ref{plot_2D_ex4} compares the exact solution and the numerical solution  computed by  our ParaDIAG-ROM preconditioner at the final time, where the boundary layer was accurately approximated without any obvious spurious oscillations.
Due to the boundary layer, our ParaDIAG-ROM preconditioner requires even larger dimension $\bar r$ (and slightly more V-cycles) for accurate approximation, which still deliver much better computational efficiency than the ParaDIAG-MG preconditioner.

\begin{table}[!htb]
	\centering 
	\caption{Results of preconditioned FGMRES for Example 2d: 2D C-D equation with upwind difference ($\epsilon=1/200, T=10$)}
	\begin{tabular}{|c||ccc||cccc|ccc|cccc|}\hline
		&\multicolumn{3}{c||}{ParaDIAG-MG(ILU) Preconditioner}  &\multicolumn{4}{c|}{ParaDIAG-ROM Preconditioner}
		\\
		\hline
		$(N,N,K)$& Rel. Res.  & Iter &CPU& Rel. Res.   & Iter  &CPU & $\bar r$ (Vc)\\   \hline 
($64^2$,640)&		 2.1E-07&	   9&	 13.5(11.9)&	    1.2E-08&	  4&	 4.1 (3.4)&	34 (7) \\
($128^2$,1280)&		 2.8E-07&	   10&	 116.9(99.0)&	    1.1E-08&	  4&	 26.4 (20.3)&	41 (9) \\
($256^2$,2560)&		 6.0E-07&	   10&	 900.6(742.4)&	    3.7E-07&	 	 3&	 124.4 (91.7)&	42 (10) \\
		\hline
	\end{tabular}
	\label{Tab_2D_ex4}
\end{table}

 \section{Conclusion}\label{finalSec}
In this work, we developed a ROM-accelerated parallel-in-time preconditioner for solving all-at-once systems from evolutionary PDEs. 
The ROM is used to reduce the computational cost for solving the sequence of complex-shifted linear systems arising from Step-(b) of the ParaDIAG algorithm. 
The reduced basis method is applied to build the ROM online, in which a greedy basis selection algorithm with several algorithmic improvements is used for finding the reduced basis efficiently. 
A variety of numerical examples are tested that illustrate the efficacy of the proposed ROM-accelerated ParaDIAG preconditioner. 
Compared with the state-of-the-art multigrid-based ParaDIAG preconditioner, the proposed approach gains more than an order of magnitude speed-up in CPU times, although the former formally has the ``optimal" complexity.
The dimension of the reduced basis seems to be moderately problem dependent,
where the convection-dominated cases require higher ROM dimensions than the diffusion-dominated ones.
Numerical investigations on the performance and scalability of our proposed ParaDIAG-ROM preconditioner in parallel computing for 3D problems are currently undertaken and will be reported elsewhere.

We plan to investigate more efficient RBM methods such as the reduced collocation methods \cite{chen2013reduced} at the next step.  
We have mainly focused on the backward Euler scheme in time, which has the first-order accuracy, but the generalization of our proposed approach to the second-order or higher-order time schemes (i.e., BDF2 scheme \cite{wu2020parallel}) is straightforward.
 The application of our current approach to the more difficult hyperbolic wave equation \cite{gander2019direct,WL2020SIMAX,de2020convergence} and its optimal control problem \cite{WL2020SISC} is another ongoing work,
 where the indefiniteness of complex-shifted linear systems in Step-(b) leads to dramatically enlarged reduced basis dimensions and hence results in less efficiency.
 It is also valuable to generalize our approach to treat similar complex-shifted linear systems with varying right-hand-sides arising in the  Laplace transform-based parallelizable contour integral method (see e.g., \cite{sheen2000parallel,sheen2003parallel,in2011contour,pang2016fast,guglielmi2020numerical}).

 \section*{Acknowledgements}
 Z.W. gratefully acknowledges support from U.S. National Science Foundation through grant DMS-1913073.
  
\bibliographystyle{spmpsci}      
	\bibliography{waveControl,Sinc,rom} 
	
\end{document}